\documentclass[11pt]{amsart}

\usepackage{hyperref}
\usepackage{amssymb,amsmath}
\usepackage[alphabetic,nobysame,abbrev]{amsrefs}
\usepackage{enumerate}
\usepackage{color}
\usepackage{pdfsync}
\RequirePackage{fix-cm}
\usepackage{bm}
\usepackage{mdframed}
\usepackage{cancel}
\usepackage{makecell}
\usepackage{pdfpages}
\usepackage{pgf,tikz}
\usepackage{mathrsfs}
\usetikzlibrary{arrows}
\usetikzlibrary{arrows,decorations.markings}
\usetikzlibrary{matrix}
\usetikzlibrary{backgrounds}

\definecolor{qqqqff}{rgb}{0.,0.,1.}
\definecolor{xdxdff}{rgb}{0.49019607843137253,0.49019607843137253,1.}

\definecolor{qqqqff}{rgb}{0.,0.,1.}

\usepackage{tikz}
\usepackage{graphicx}
\usepackage[cmtip,all]{xy}

\usepackage[draft]{todonotes}
\usepackage[cmtip,all]{xy}

\usepackage{hyperref}
\usepackage[utf8]{inputenc}
\usepackage{amsmath,amsfonts,amssymb,euscript,graphicx,color,tikz}

\theoremstyle{plain}
\newtheorem{theorem}{Theorem}[subsection]
\newtheorem{thm}[theorem]{Theorem}
\newtheorem{lem}[theorem]{Lemma}
\newtheorem{cor}[theorem]{Corollary}
\newtheorem{pro}[theorem]{Proposition}

\theoremstyle{definition}
\newtheorem{DEF}[theorem]{Definition}

\newtheorem{exa}[theorem]{Example}
\newtheorem{con}[theorem]{Convention}

\newtheorem{rem}[theorem]{Remark}

\newtheorem{parag}[theorem]{{}}

\newcommand\pref[1]{\textbf{\ref{#1}}}
\numberwithin{equation}{section}

\newcommand{\sub}{\subseteq}

\newcommand{\la}{Lie algebra }

\newcommand{\fm}{(\cdot,\cdot)}

\newcommand{\ep}{\hfill$\Box$}

\def\ad{\hbox{ad}}
\def\andd{\quad\hbox{and}\quad}
\def\sg{\sigma}
\def\a{\alpha}
\def\b{\beta}
\def\lam{\lambda}
\def\Lam{\Lambda}
\def\ep{\epsilon}
\def\andd{\quad\hbox{and}\quad}

\def\id{\hbox{id}}
\def\Aut{\hbox{Aut}}

\def\andd{\quad\hbox{and}\quad}
\def\ind{\hbox{ind}}

\def\v{{\mathcal V}}
\def\u{{\mathcal U}}
\def\vd{\dot{\mathcal V}}

\def\vt{\tilde{\mathcal V}}

\def\fm{(\cdot,\cdot)}
\def\a{\alpha}

\def\w{{\mathcal W}}

\def\sub{\subseteq}
\def\rd{\dot{R}}

\def\lam{\lambda}
\def\Lam{\Lambda}

\def\1k{\frac{1}{k}}

\def\la{\langle}
\def\ra{\rangle}
\def\rds{\dot{R}_{sh}}
\def\rdl{\dot{R}_{lg}}

\def\GL{GL}

\def\d{\delta}

\def\b{\beta}

\def\qed{\hfill$\Box$}

\def\sg{\sigma}

\def\hh{{\mathcal H}}

\def\sg{\sigma}

\usepackage{verbatim} 

\def\ad{\hbox{ad}}

\def\bbbc{{\mathbb C}}
\def\bbbz{{\mathbb Z}}

\def\bbbr{{\mathbb R}}

\def\bbbk{{\mathbb K}}

\def\cc{{\mathcal C}}

\def\ep{\epsilon}

\def\ii{{\mathcal I}}
\def\proof{{\noindent\bf Proof. }}

\def\rds{\dot{R}_{sh}}
\def\rdl{\dot{R}_{lg}}

\def\span{\hbox{span}}

\def\pp{\mathcal P}

\def\refl{\hbox{refl}}
\def\rank{\hbox{rank}}
\def\null{\hbox{null}}

\def\DynkinNodeSize{1.5mm}
\def\DynkinArrowLength{2mm}
\tikzset{
dnode/.style={
circle,
inner sep=0pt,
minimum size=\DynkinNodeSize,
fill=white,
draw},
middlearrow/.style={
decoration={markings,
	mark=at position 0.8 with
	{\draw (0:0mm) -- +(+140:\DynkinArrowLength); \draw (0:0mm) -- +(-140:\DynkinArrowLength);},
},
postaction={decorate}
},
leftrightarrow/.style={
decoration={markings,
	mark=at position 0.999 with
	{
		\draw (0:0mm) -- +(+135:\DynkinArrowLength); \draw (0:0mm) -- +(-135:\DynkinArrowLength);
	},
	mark=at position 0.001 with
	{
		\draw (0:0mm) -- +(+45:\DynkinArrowLength); \draw (0:0mm) -- +(-45:\DynkinArrowLength);
	},
},
postaction={decorate}
},
sedge/.style={
},
dedge/.style={
middlearrow,
double distance=0.6mm,
},
tedge/.style={
middlearrow,
double distance=1.0mm+\pgflinewidth,
postaction={draw}, 
},
infedge/.style={
leftrightarrow,
double distance=0.5mm,
},
}

\begin{document}

%
%
\title{Characters for extended affine Lie algebras; a combinatorial approach}

\author{Saeid Azam}
\address
{Department of Pure Mathematics\\Faculty of Mathematics and Statistics\\
	University of Isfahan\\ P.O.Box: 81746-73441\\ Isfahan, Iran, and\\
	School of Mathematics, Institute for
	Research in Fundamental Sciences (IPM), P.O. Box: 19395-5746.} \email{azam@ipm.ir, azam@sci.ui.ac.ir}
\thanks{This work is based upon research funded by Iran National Science Foundation (INSF) under project No.  4001480.}
\thanks{This research was in part carried out in
	IPM-Isfahan Branch.}
\keywords{\em Extended affine Lie algebra, Character, Cartan automorphism, Chavalley involution, involutive Lie algebra}

\begin{abstract}
	The behavior of objects associated with general extended affine Lie algebras is typically distinct from their counterparts in affine Lie algebras. Our research focuses on studying characters and Cartan automorphisms, which appear in the study of Chevalley involutions and Chevalley bases for extended affine Lie algebras. We show that for almost all extended affine Lie algebras, 
	any finite order Cartan automorphism is diagonal, and its corresponding combinatorial map is a character. 
	\end{abstract}
 \subjclass[2020]{17B67, 17B22, 16W10, 17B65}
\maketitle

\section{Introduction}\setcounter{equation}{0}\label{preliminaries}
\markboth{S. Azam}{Charatcers for extended affine Lie algebars}

Extended affine Lie algebras are a class of Lie algebras that are defined by a set of axioms. They share many common features with finite and affine Lie algebras, but also possess unique characteristics that make them fascinating objects of study, see \cite{H-KT90}, \cite{AABGP97} and \cite{Neh11}. For instance, like finite and affine Lie algebras, an extended affine Lie algebra is equipped with a symmetric invariant non-degenerate bilinear form, and has a finite-dimensional Cartan subalgebra that splits the Lie algebra into a direct sum of root spaces.  

The set of roots of an extended affine Lie algebra is called an extended affine root system. Similar to affine Lie algebras, a root is called isotropic if any integral multiple of it is also a root. The dimension of the vector space spanned by isotropic roots is known as the nullity of an extended affine Lie algebra or its root system. 
The finite and affine Lie algebras correspond to nullities 0 and 1, respectively.

Concerning the structure theory, two of the most important objects that are associated with an extended affine Lie algebra are the core which is the subalgebra generated by non-isotropic root spaces, and the centerless core which is the quotient algebra of the core over its center. Similarly, associated with an extended affine root system $R$ there are some objects which are essential in its study, we mention some here which are more related to our work.
The image $\bar R$ of $R$ in the quotient space of the real span of roots over the span of isotropic roots is an irreducible finite root system whose rank and type are turned to be isomorphism invariants of $R$.  These are called the rank and the nullity of $R$, respectively. In this work, we assume that $R$ is reduced. The root system $R$ can be then described with three ingredients; an {irreducible finite root system} $\rd$ isomorphic to $\bar R$, and two discrete subsets $S$ and $L$ of isotropic roots. The interrelations between $S$ and $L$, called semilattices, reflect the interaction of reflections based on short and long roots of $R$.

In this study, we investigate the concepts of "characters" and ``diagonal automorphisms" for extended affine Lie algebras. These concepts exhibit characteristics in higher nullities {which} are absent in finite or affine cases. Our study of characters stems from our interest in Chevalley involutions. These are finite order automorphisms that act as minus identity on the Cartan subalgebra.
Chevalley involutions are essential in the study of Chevalley bases, and more generally in the modular theory of extended affine Lie algebras, see \cite{AFI22} and \cite{AI23}. The Composition of two Chevalley involutions results in an automorphism that stabilizes pointwise the Cartan subalgebra, which we refer to as a Cartan automorphism. Involutions, in a more broad context have been studied in the context of involutive Lie algebras, see for example \cite{Gao96} and \cite{Cal04}. 

Let $E$ be an extended affine Lie algebra with Cartan subalgebra $\hh$ and root system $R$.
Since a Cartan automorphism $\Psi$ stabilizes each root space, it acts on the root space $E_\a$, $\a$ non-isotropic, as $\psi(\a)\id$ where $\psi(\a)$ is a non-zero complex number. If $\Psi$ is diagonal, i.e., if it also acts as a scalar multiple of identity on each isotropic root space, then the map   
$\psi:R\rightarrow\bbbc^\times$ is called a core-character associated with $\Psi$ for $R$.
If $\psi(\a+\b)=\psi(\a)\psi(\b)$ for $\a,\b,\a+\b\in R$, then $\psi$ is called a character, and if $\psi$ is extendable to the root lattice $\la R\ra$, we say that $\psi$ is a root-lattice character. While any character for $R$ leads to a diagonal automorphism, the converse is not necessarily true. In this work, we investigate the interrelations between Chevalley involutions, Cartan automorphisms, and characters. The paper is organized as follows.

In Section 2,  we provide the definition of an extended affine Lie algebra denoted by $E$ and report some fundamental facts about $E$ and its root system. In particular, we explain the concepts of ``index'' and ``reflectable base''  which are essential to extended affine root systems.

{In Section \ref{sec9}, we establish certain basic facts about Chevalley involutions.} Specially, we prove that any finite order Chevalley involution has order $2$, as stated in Lemma \ref{l1}, and that a finite-order automorphism whose restriction to the core is a Chevalley involution for the core is also a Chevalley involution for $E$ only if it stabilizes the Cartan subalgebra, see Corollary \ref{chev2}. Furthermore, we show that the existence of Chevalley involutions does not depend on the choice of the Cartan subalgebra, see Lemma \ref{cong}. 

In Section 4, the core part of the paper, we consider an extended affine Lie algebra $E$ of type $X$ and $\rank>1$, with root system $R$; also consider a finite order Cartan automorphism $\Psi$ of $E$ and let $\psi:R\rightarrow\bbbk^\times$ be the core-character associated with $R$, see Definition \ref{aug31a}. We proceed with a type-dependent combinatorial approach to show that if $X\not = B_\ell$, or if $X=B_\ell$ of index $0$, then $\Psi$ is diagonal and $\psi$ is a character, see Propositions \ref{lemaug31a}, \ref{bim11}, \ref{bim15}, and \ref{aug29a}.
The connections among Chevalley involutions, Cartan automorphisms, and characters are summarized in Theorem \ref{cor021}.

In Section 5, the final section, we show that in contrast to the affine case, not every character for an extended affine root system is expendable to a root-lattice character, see Definition \ref{bims31} and  Corollary \ref{diag10}.

\section{\bf Preliminaries}\setcounter{equation}{0}\label{preliminaries}
Throughout this work $\bbbk$ is the field of complex numbers, and all vector spaces are over $\bbbk$ unless otherwise is specified. We set $\bbbk^\times=\bbbk\setminus\{0\}$.
We denote the dual space of a vector space $\u$ with $\u^\star$.  The notation $R_1\uplus R_2$ means 
the union of two disjoint sets $R_1$ and $R_2$. We denote by $\la S\ra$, the $\bbbz$-span of a subset $S$ of a vector space.

\subsection{Extended affine Lie algebras and root systems}

We begin by recalling some terminologies and basic facts about extended affine Lie algebras and extended affine root systems. For more details, we refer the reader to \cite{AABGP97} and \cite{Neh11}. 

\begin{parag}\label{july18} An {\it extended affine Lie algebra} is a triple $(E,\fm,\hh)$ where
	$E$ is a Lie algebra over $\bbbk$, $\hh$ is a non-trivial finite dimensional subalgebra of $E$, and $\fm$ is a $\bbbk$-valued symmetric bilinear form on $E$ satisfying axioms (A1)-(A6) below:
	
	(A1) The bilinear form $\fm$ is invariant and non-degenerate.
	
	(A2) $\hh$ is a  splitting Cartan subalgebra of $E$ meaning that  
	$E=\sum_{\a\in R}E_\a$ with
	$E_\a=\{x\in E\mid [h,x]=\a(h)x\hbox{ for all }h\in\hh\},$
	$R=\{\a\in\hh^\star\mid E_\a\not=\{0\}\}$,
	and $E_0=\hh.$ The set $R$ is called the {\it root system} of $E$.

	From axioms (A1)-(A2), we see that the form $\fm$ on $E$ restricted to $\hh$ is non-degenerate. For $\a\in\hh^\star$ let $t_\a\in\hh$
	be the unique element satisfying $\a(h)=(h,t_\a)$, $h\in\hh$, and consider the form on $\hh^\star$ defined by $(\a,\b):=(t_\a,t_\b)$.
	The root $\a$ is called {\it isotropic} if
	$(\a,\a)=0$ and {\it non-isotropic} otherwise.
	We denote by $R^0$ the set of isotropic roots and by $R^\times$ the set of non-isotropic roots of $R$. The subalgebra $E_c$ of $E$ generated by non-isotropic root spaces is called the {\it core} of $E$.
	
	(A3) For $\a\in R^\times$ and $x\in E_\a$, $\ad(x)$ is locally nilpotent on $E$.
	
	(A4) $E$ is {\it tame}, that is, the centralizer of $E_c$ in $E$ is contained in $E_c$.
	
	(A5) The $\bbbz$-span of $R$ in $\hh^\star$ is a free abelian group of finite rank.
	
	(A6) $R^\times$ is {\it indecomposable} meaning that if  $R^\times=R_1\uplus R_2$ with $R_1\not=\emptyset$ and $R_2\not=\emptyset,$ then  $R_1\not\perp R_2$. 
	%
\end{parag}


\begin{rem}\label{remkhan1}
	In order to better achieve our goals, we have chosen to include the tameness condition (axiom (A4)) in the definition of an extended affine Lie algebra, despite it is typically not included in the definition. It follows from the tameness condition that isotropic roots are non-isolated, meaning that for $\sg\in R^0$, there exists $\a\in R^\times$ with $\a+\sg\in R$.
\end{rem}

	Assume that $(E,\fm,\hh)$ is an extended affine Lie algebra.
	The term
	$$h_\a:=\frac{2t_\a}{(\a,\a)}\qquad(\a\in R^\times),$$ will be used frequently in the sequel.
	It follows that for $\a\in R^\times$, there exists $x_{\pm\a}\in E_{\pm\a}$ such that $(x_\a,h_\a,x_{-\a})$ forms an $\frak{sl}_2$-triple.
	Also we recall from \cite[Remark 1.5(ii)]{Az06} and \cite[Chapter I]{AABGP97} that
	\begin{equation}\label{Aug21a}
	[E_\a,E_\b]\not=\{0\}\qquad(\a\in R^\times,\;\b,\a+\b\in R),
	\end{equation}
\begin{equation}\label{city}
	[x_\a,x_{-\a}]=(x_\a,x_{-\a})t_\a\qquad(\a\in R, {x_{\pm\a}\in E_{\pm\a}}).
\end{equation}
	
	\begin{parag}\label{Aug20a}
	Next, we recall some facts about the root system $R$ of $E$ which is called an {\it extended affine root system}. Set $\v:=\span_{\bbbr}R$ and $\v^0:=\span_{\bbbr}R^0$. 
		It turns out that the form on $\hh^\star$ restricted to $\v$ can be assumed to be positive semidefinite, $2\a\not\in R^\times$ for $\a\in R^\times$, and that the root-string property holds, namely for $\a\in R^\times$ and $\b\in R$, there exist non-negative integers $d,u$ such that
		\begin{equation}\label{string}
			\{n\a+\b\in R\mid n\in\bbbz\}=\{\b-d\a,\ldots,\b,\ldots\b+u\a\}
			\end{equation} with
		$2(\b,\a^\vee)/(\a,\a)=d-u.$ Moreover, 
		the image $\bar R$ of $R$ in $\bar\v:=\v/\v^0$ is an irreducible finite root system in $\bar\v$. One can find a preimage $\dot R$ of $\bar R$ in $\v$, under $\bar{\;}$, such that $\dot R$ is an irreducible finite root system in its real span, isomorphic to $\bar R$. Then $R\sub\dot R+\Lam$, where $\Lam=\la R^0\ra$ is a free abelian group of rank equal to $\dim\v^0$.

	The {\it rank} and the {\it type} of $E$ or $R$ is by definition the rank and the type of $\rd$. The dimension of $\v^0$ is called the {\it nullity} of $R$. We denote the rank and the nullity of $R$ by $\hbox{rank}(R)$ and $\null(R)$, respectively.
	We note that $\rank \la R\ra=\rank(R)+\null(R)$.
	In this work, {\it we assume that $E$ is of reduced type}.
\end{parag}

\subsection{Semilattices}	
	The internal structure of an extended affine root system $R$ can be described in terms the finite root system $\rd$ and certain subsets of $R^0$ called "semilattices". 
	
	\begin{parag}\label{bims26}
		A {\it semilattice} in a $\nu$-dimensional real vector space $\u$ is a discrete spanning subset $S$ containing $0$ with $S\pm2S\sub S$. It turns out that $\Lam:=\la S\ra$ is a lattice in $\u$, and $S$ is the union of a set of cosets of $2\Lam$ in $\Lam$ including the trivial coset $2\Lam$. Moreover,  $S=\uplus_{i=0}^m(\tau_i+2\Lam)$, where $\tau_i$'s can be chosen such that, $\tau_i$'s represent distinct cosets of $2\Lam$, $\tau_0=0$, and $\Lam=\sum_{i=1}^\nu\bbbz\tau_i$. 
		The integer $m$,
	the number of non-trivial cosets of $2\Lam$ in $S$, is called the {\it index} of $S$, denoted by $\ind(S)$. We call $\nu$ the {\it rank} of $S$. Note that  $\nu\leq \ind(S)\leq 2^\nu$.  
	\end{parag}

\begin{parag}\label{bims23}
		Now, we can give the following description of $R$, in terms of the length of roots and two involved semilattices, 
	\begin{equation}\label{bims1}
	\begin{array}{c}
	R^0=S+S\andd R=R^0\cup (\rds +S)\cup (\rdl+L),\\
	S+L=S,\quad kS+L=L,
	\end{array}
	\end{equation}
	where $\rds$ and $\rdl$ are the sets of short and long roots of $\rd$, respectively, and $k$ is the maximum number of edges appearing between two consecutive roots in the Dynkin diagram of $\rd$. Here if $\rdl=\emptyset$, we ignore in (\ref{bims1}) all terms in which $\rdl$ or $L$ appear. It turns out that
	
	- If $R$ is simply laced of rank $>1$, or is of type $F_4$ or $G_2$, then $S$ and $L$ are lattices.
	
	- If $R$ is of type $B_\ell$ (respectively $C_\ell$), $\ell\geq 3$, the $L$ (respectively $S$) is a lattice.
	
	To have a more precise description of $R$ when it is of one of non-simply laced types, we need to introduce another isomorphism invariant of $R$ called the {\it twist number}, given by $k^t=|\la L\ra/\la S\ra|$. If $R$ has twist number $t$, then there exists semilattices $S_1$ of rank $t$ and $S_2$ of rank $\nu-t$, such that
	\begin{equation}\label{bims23}
		S=S_1\oplus\la S_2\ra\andd L=k\la S_1\ra\oplus S_2.
		\end{equation}
		\end{parag}


		Next, we explain briefly about the Weyl group of $R$. Set $\vt:=\v\oplus(\v^0)^\star$, where $(\v^0)^\star$ is the dual space of $\v^0$. The form $\fm$  on $\v$ extends to a non-degenerate form on $\vt$ by 
	$$ (\dot\v,(\v^0)^\star)=((\v^0)^\star,(\v^0)^\star)=\{0\},\;
	(\gamma,\sg):=\gamma(\sg),\quad 
	(\gamma\in(\v^0)^\star, \sg\in\v^0).$$
	For $\a\in R^\times$ the reflection $w_\a\in\GL(\vt)$ is defined by
	$w_\a(\b)=\b-(2(\a,\b)/(\a,\a))\a.$ The {\it Weyl group} $\w$ of $R$ is the subgroup of $\GL(\vt)$ generated by reflections $w_\a$, $\a\in R^\times.$

\begin{parag}\label{july20b}
For $\pp\sub R^\times$, denote by $\w_\pp$ the subgroup of $\w$ generated by reflections $w_\a$, $\a\in\pp$. A {\it reflectable base} for $R$ is a subset $\pp$ of $R^\times$ such that $\w_\pp\pp=R^\times$ and no proper subset of $\pp$ has this property. According to \cite[Corollary 4.2.2]{APT23}, any reflectable base is finite. 
	
	We define $\refl(R)$ to be the smallest cardinality of a reflectable base for $R$, and $\ind(R)$ as $\refl(R)-\rank(R)-\null(R)=\refl(R)-\rank\la R\ra$. These terms, along with $\rank(R)$ and $\null(R)$, are all isomorphism invariants for $R$.
	We may also describe the $\ind(R)$ in terms of rank $\ell$, nullity $\nu$, and semilattices $S,$  $S_1$ and $S_2$ given in (\ref{bims1}) and (\ref{bims23}). In fact, depending on type $X$ of $R$, we have
	\begin{equation}\label{bims25}
		\ind(R)=\left\{\begin{array}{ll}
			\ind(S)-\ell-\nu&X= A_1,\\
			0&X=A_\ell\;(\ell\geq 2),D_\ell, E_{6,7,8}, F_4, G_2,\\
			\ind(S_1)+\ind(S_2)-\nu& X=B_2,\\
			\ind(S_1)-t&X=B_\ell,\;\ell\geq 3,\\
			\ind(S_2)-\nu-t&X=C_\ell,\;\ell\geq 3.
		\end{array}\right.
	\end{equation}
	Note that in type $B_\ell$, $\ell\geq 3$, if $\ind(R)=0$, then $\ind(S_1)=t=\rank (S_1)$, so as we saw in \pref{bims26}, in this case we have
	$S_1=\cup_{i=0}^t(\tau_i+2\la S_1\ra)$ with $\la S_1\ra=\sum_{i=1}^t\bbbz\tau_i$.
	
	Extended affine root systems of index $0$ are more similar to finite and affine root systems than those with a non-zero index, as per \cite{APT23}. We note from \cite[\S 4]{Az99} that $\ind(R)=0$ if $R$ is finite or affine (i.e., if $\null(R)=0$ or $1$) or if $R$ has one of the types $A_\ell$ ($\ell\geq 2$), $D_\ell$ ($\ell\geq 4$), $F_4$, $G_2$, $E_{6,7,8}$. Therefore, the extended affine root systems with a non-zero index occur only in types $A_1$, $B_\ell$, and $C_\ell$.
\end{parag}


\section{\bf Chevalley involutions}\setcounter{equation}{0}\label{sec9}
We begin the section by establishing some fundamental facts regarding Chevalley involutions for an extended affine Lie algebra $E=(E,\fm,\hh)$.
We set $E_{c,\a}:=E_c\cap E_\a$ for $\a\in R$. When we refer to an involution on a vector space, we are referring to an automorphism of order $2$.



\subsection{Chevalley involutions}
Although we have assumed that all extended affine Lie algebras in this work are reduced, the results in this subsection are also valid for non-reduced ones.

\begin{lem}\label{l2}
Suppose $\tau$ is a finite order automorphism of an extended affine Lie algebra. Then $\tau=\id$ if and only the restriction of $\tau$ to the core is identity.
\end{lem}

\proof Suppose $\tau_{|_{E_c}}=\id_{E_c}$. Let $x\in E$ and  $y\in E_c$. 
Since $E_c$ is an ideal of $E$, we get $[x,y]=\tau [x,y]=[\tau(x), y]$.
Therefore, $\tau(x)-x$ is in centralizer of the core. Since $E$ is tame, $\tau(x)-x\in E_c$ and so
$\tau(\tau(x)-x)=\tau(x)-x,$ for each $x\in E$. Thus $\tau$ satisfies the polynomial $(t-1)^2=0$. Since 
$\tau$ is of finite order, we get $\tau=\id$.\qed
\begin{lem}\label{l1}
Let $\tau$ be a finite order automorphism of $E$. Fix $\ep\in\{1,-1\}$.

(i) If $\tau(E_{\a})=E_{\ep\a}$ for all $\a\in R^\times$, then 

- $\tau_{|_{\hh\cap E_c}}=\ep\id,$

- if $\ep=-1$ then $\tau^2=\id$,

- $\tau_{|_\hh}=\ep\id_{\hh}$ mod $Z(E_c)$. 

(ii) $\tau(E_\a)=E_{\ep\a}$ for all $\a\in R$ if and only if
$\tau(h)=\ep h$ for all $h\in \hh$.

(iii) If {$\tau(E_{\a})=E_{\ep\a}$ for all $\a\in R$, then $\tau$ preserves the form.}
\end{lem}

\proof
(i) Assume that $\tau(E_\a)=E_{\ep\a}$ for each $\a\in R^\times$. For $\a\in R^\times$, fix
$x_{\pm\a}\in E_{\pm\a}$ such that $[x_\a,x_{-\a}]=h_\a$. Then
$\tau(x_\a)=\mu_\a x_{\ep\a}$ for some $\mu_\a\in\bbbk^\times$. So
$$\tau(h_\a)=\tau[x_\a,x_{-\a}]=\mu_\a\mu_{-\a}[x_{\ep\a},x_{-\ep\a}]=\mu_\a\mu_{-\a}h_{\ep\a}=
\ep\mu_\a\mu_{-\a}h_\a.$$
Now by applying $\tau$ to the both sides of $2x_\a=[h_\a,x_\a]$, we get
$$
2\mu_\a x_{\ep\a}=\ep\mu_\a\mu_{-\a}[h_\a,\mu_\a x_{\ep\a}]=\ep^2 2\mu_\a^2\mu_{-\a}x_{\ep\a}=2\mu_\a^2\mu_{-\a}x_{\ep\a}.$$ 
Thus 
\begin{equation}\label{e2}
\mu_\a\mu_{-\a}=1\andd\tau(h_\a)=\ep h_\a, \hbox{ for each }\a\in R^\times.
\end{equation}
Since $\hh\cap E_c=\sum_{\a\in R^\times}h_\a$, we get
\begin{equation}\label{e1}
\tau_{|_{\hh\cap E_c}}=\ep\id.
\end{equation}
Note that if $\ep=-1$ then form  (\ref{e2}) for $\a\in R^\times$ we have
$\tau^2(x_\a)=\mu_\a\mu_{-\a}x_\a=x_\a$.
Since $E_c$ is generated by $x_\a$, $\a\in R^\times$, we get $\tau^2=\id_{|_{E_c}}$, and so by  Lemma \ref{l2}, $\tau^2=\id$.

Next for $h\in\hh$ and $\a\in R^\times$, we have
\begin{eqnarray*}
[h-\ep\tau(h),x_\a]&=&[h,x_\a]-\ep\tau[h,\mu_{\ep\a}^{-1}x_{\ep\a}]\\
&=&\a(h)x_\a-\ep\ep\mu_{\ep\a}^{-1}\a(h)\tau(x_{\ep\a})\\
&=&\a(h)x_\a-\a(h)\mu_{\ep\a}^{-1}\mu_{\ep\a}x_\a=0.
\end{eqnarray*}
Therefore, $h-\ep\tau(h)$ is in the centralizer of $E_c$ in $E$ for each $h\in\hh$. Since $E$ is tame,
$h-\ep\tau(h)\in Z(E_c)$.

(ii) Assume first that  $\tau(E_\a)=E_{\ep\a}$ for all $\a$; in particular,
$\tau(\hh)=\tau(E_0)=E_0=\hh$. So from part (i), we have $h-\ep\tau(h)\in\hh\cap Z(E_c)$. Then from (\ref{e1}) we have
$\tau(h-\ep\tau(h))=\ep(h-\ep\tau(h))$ for each $h\in\hh$. Thus $\ep\tau$ restricted to $\hh$ satisfies the equation $(t-1)^2=0$. Since $\tau$ is of finite order, it follows that $\tau_{|_\hh}=\ep\id$.

Conversely, assume that $\tau_{|_\hh}=\ep\id.$ Then
for $\a\in R$, $h\in\hh$ and $x\in E_\a$, we have 
$$[h,\tau(x)]=\tau[\ep h,x]=\ep\a(h)\tau(x),$$
and so $\tau(x)\in E_{\ep\a}.$

(iii) By (ii), {for $h,h'\in\hh=E_0$}, we have
$(\tau(h),\tau(h'))=(\ep h,\ep h')=(h,h')$. Assume now that $0\not=\a\in R$, $\b\in R$.
Let {$x_\a\in E_\a,x_\b\in E_\b$}. If $\a+\b\not=0$, then
$(\tau(x_\a),\tau(x_\b))=0=(x_\a,x_\b)$. If $\b=-\a$, then by (ii) and (\ref{city}),
	$$t_{\ep\a}(\tau(x_\a),\tau(x_{-\a}))=\tau[x_\a,x_{-\a}]=t_{\ep\a}(x_\a,x_{-\a}).$$
Therefore, $(x_\a,x_\b)=(\tau(x_\a),\tau(x_\b))$ for $\a,\b\in R$ and {$x_\a\in E_\a,x_\b\in E_\b$.}
\qed

\begin{DEF}\label{chev} 
(i) We call a finite order automorphism  $\tau$ of $E$ a {\it Chevalley involution} if
$\tau(E_\a)=E_{-\a}$ for all $\a\in R$, or equivalently (see \ref{l1}(ii)) if $\tau(h)=-h$ for all $h\in \hh$.

(ii) We call a finite order automorphism $\tau$ of $E_c$ a {\it Chevalley involution} for $E_c$ if  
$\tau(E_\a)=E_{-\a}$ for all $\a\in R^\times$.
\end{DEF}

\begin{rem}\label{Aug12a}
Since by \cite[Proposition 3.4]{CNPY16} any automorphism $\tau$ for $E$ stabilizes the core, it is immediate from Definition \ref{chev} that any Chevalley involution for $E$ restricts to a Chevalley involution for
$E_c$. Also, we have from Lemma \ref{l1}(iii) that any finite order Chevalley involution preserves the form.
\end{rem}

\begin{cor}\label{chev2}
Let $\tau$ be a finite order automorphism of $E$ such that the restriction of $\tau$ to $E_c$ is a Chevalley involution for $E_c$. Then
$\tau$ is a Chevalley involution for $E$ if and only if $\tau$ stabilizes $\hh$.
In particular, if $E$ is an affine Kac-Moody Lie algebra, then $\tau$ is a Chevalley involution for $E$ if and only if its restriction to $E_c$ is a Chevalley involution. 
\end{cor} 

\proof
Since $\tau$ is a Chevalley involution for $E_c$, we have $\tau(E_\a)=E_{-\a}$ for $\a\in R^\times$. By assumption $\tau$ stabilizes $\hh=E_0$, so we have $\tau(E_\a)=E_{-\a}$ for all $\a\in R$. Thus the first part of the statement is an immediate consequence of Lemma \ref{l1}(ii).
For the second part, we note that if $E$ is an affine Kac--Moody Lie algebra, one knows that $Z(E_c)\sub\hh$, see \cite{Kac90}, thus $\tau$ stabilizes $\hh$ and we are done.\qed

\begin{lem}\label{cong}
The existence of Chevalley involutions for an extended affine Lie algebra does not depend on the choice of the Cartan subalgebra.
\end{lem}

\proof
Suppose $\hh$ and $\hh'$ are two Cartan subalgebras for $E$, and $\tau$ is a Chevalley involution for $(E,\fm,\hh)$. By conjugacy theorem for extended affine Lie algebras (see \cite[Theorem 7.6]{CNPY16} and
\cite{CNP17}), there exists an automorphism $\sg$ of $E$ which maps $\hh$ onto $\hh'$.
Then $\bar\tau:=\sg\tau\sg^{-1}$ is a Chevalley involution for $(E,\fm,\hh')$.
\qed

\section{Cartan automorphisms, and characters}\setcounter{equation}{0}\label{aug30b}	

Assume that $E=(E,\fm,\hh)$ is an extended affine Lie algebra with root system $R$. The composition of two Chevalley involutions for $E$ is referred to as a "Cartan automorphism" since it stabilizes pointwise the elements of the Cartan subalgebra $\hh$ and thus, stabilizes each root space. Moreover, the composition of a Cartan automorphism and a Chevalley involution results in another Chevalley involution. In turn, Cartan automorphisms can be constructed using the "characters" of the root system, and vice versa, to each Cartan automorphism one associates a map that is perceivable to be generally a character. This section studies the relationships between Chevalley involutions, Cartan automorphisms, and characters. 

\subsection{\bf Characters}
We begin by giving definitions of Cartan automorphisms, characters, and some related terms.
\begin{DEF}\label{diag} Let $(E,\fm,\hh)$ be an extended affine Lie algebra and $R$ be an extended affine root system.
	
(i) A Cartan automorphism $\psi$ for $E$ is the one that satisfies $\psi(h)=h$ for all $h\in\hh$. It is evident that Cartan automorphisms stabilize each root space.

(ii) We call a Cartan automorphism {\it diagonal} if it acts as a scalar multiple of identity when restricted to each root space.

%
	\end{DEF}

\begin{DEF}\label{bims31}  
	(i) We call a map $\psi:R\rightarrow\bbbk^\times$, a {\it core-character} for $R$ if
	$\psi(\a+\b)=\psi(\a)\psi(\b)$, for $\a\in R^\times,$ $\b,\a+\b\in R$.
	A core-character $\psi$ is called of {\it finite-order} if there exists $m\in\bbbz_{>0}$ with $\psi(\a)^m=1$ for all $\a\in R$. 
	If $R$ has type $X$, we call a core-character for $R$ a {\it core-character of type} $X$.
	
	(ii) A core-character $\psi$ is called a {\it character} if
	$\psi(\a+\b)=\psi(\a)\psi(\b)$ for $\a,\b,\a+\b\in R$.  
	A {\it root-lattice character} for $R$ is a group homomorphism $\psi:\la R\ra\rightarrow\bbbk^\times$. 
\end{DEF}

{\begin{exa}\label{sky4}
		(i)  Given an extended affine Lie algebra $E$ with root system $R$ and a (finite-)order character $\Phi:R\rightarrow\bbbk^\times$, the linear map defined by $\bar\Phi(x)=\Phi(\a)x$, for $x\in E_\a$, $\a\in R$, defines a (finite order) diagonal automorphism for $E$. 
		
		(ii) Let $m$ be a positive integer. Given an extended affine root system $R$, fix a basis $\pp$ for the free abelian group $\la R\ra$. For $\a\in\pp$ fix an $m$-th root of unity $\eta_\a$.
		Then the group homomorphism $\Phi:\la R\ra\longrightarrow\bbbk^\times$ induced by $\Phi(\a)=\eta_\a$, $\a\in\pp$, defines a finite-order root-lattice character for $R$.
\end{exa}}

\begin{rem}\label{diag8} (i) It is immediate that a Cartan automorphism stabilizes each root space.
	Conversely, suppose $\psi$ is a finite order automorphism of $E$ that stabilizes each root space. Then by Lemma \ref{l1}(ii), $\psi$ is a Cartan automorphism.

	(ii)	In \cite{L88}, the author classifies involutive automorphisms (automorphisms of order 2) for affine Kac-Moody Lie algebras by dividing them into two categories: first and second kind automorphisms. Chevalley involutions belong to the second kind, while Cartan involutions belong to the first kind.
	
	(iii) By Lemma \ref{l1}(iii), any finite order Cartan automorphism preserves the form. 
	\end{rem}

Previously, in Example \ref{sky4}, we learned about producing non-trivial finite order diagonal automorphisms for the extended affine Lie algebra $E$ using a finite-order character for the root system $R$. But does every finite order diagonal automorphism arise this way? We will find out in the following sections that the answer is "almost" affirmative, with a few exceptions in types $A_1$ and $B_\ell$.

From now on we assume that $E$ is an extended affine Lie algebra with root system $R$, and we fix a Cartan automorphism $\Psi$. Hereafter, whenever we write $x_\a$, $\a\in R$, we mean $x_\a$ is a non-zero element of $E_\a$.

For $\a\in R^\times$, we define
$\eta_\a$ by
\begin{equation}\label{aug31b}
\Psi(x_\a)=\eta_\a x_\a,\quad (0\not=x\in E_\a,\a\in R^\times).
\end{equation}
Now let $\sg\in R^0$. Since isotropic roots are non-isolated, there exists $\a\in R^\times$ with
$\a+\sg\in R^\times$. We define
$$\eta_\sg^\a:=\eta_{\a+\sg}\eta_{-\a},\quad (\sg\in R^0,\a,\a+\sg\in R^\times).
$$

\begin{con}
	Whenever we use the notation $\eta^\a_\sg$, we understand that $\a\in R^\times$, $\sg\in R^0$, {$\a+\sg\in R^\times$.}
\end{con}

\begin{lem}\label{bims30}
	Let $\Psi$ be a finite order Cartan automorphism. Then $\Psi$ is diagonal if and only if $\eta^\a_\sg=\eta^\b_\sg$,
	$\a,\b\in R^\times$, $\sg\in R^0$.
	\end{lem}

\proof Assume first that $\Psi$ is diagonal. Then
$\Psi_{|_{E_\a}}=\eta_\a \id_{|_{E_\a}}$, $\a\in R$, for some $\eta_\a\in\bbbk^\times$.
Let $\sg\in R^0$ and $\a,\a+\sg\in R^\times$. By (\ref{Aug21a}), there exist
$x_\a\in E_\a$ and $x_\sg\in E_\sg$ such that $[x_\a,x_\sg]\not=0$. Again by
(\ref{Aug21a}),
$X_\sg:=[x_{-\a},[x_\a,x_\sg]]\not=0$. Then
\begin{eqnarray*}
	\eta_\sg [x_{-\a},[x_\a,x_\sg]]
	&=&
	\Psi([x_{-\a},[x_\a,x_\sg]])\\
	&=&[\Psi(x_{-\a}),\Psi[x_\a,x_\sg]]\\
	&=&
	\eta_{-\a}\eta_{\a+\sg}[x_{-\a},[x_\a,x_\sg]]\\
	&=&
	{\eta_\sg^\a [x_{-\a},[x_\a,x_\sg]].}
	\end{eqnarray*}
This shows that the definition of $\eta^\a_\sg$ is independent of the choice of $\a\in R^\times.$

Conversely, assume that $\eta_\sg^\a=\eta_\sg^\b$ for $\sg\in R^0$, $\a,\b\in R^\times.$ So for $\sg$ in $R^0$, the term
$\eta_\sg:=\eta_\sg^\a$, for some $\a\in R^\times$, is well-defined. 
Since for $\a\in R^\times$, $\dim E_\a=1$, we have $\Psi_{|_{E_\a}}=\eta_\a\id_{|_{E_\a}}$. Now let $\sg\in R^0$. Then if $\a,\a+\sg\in R^\times$, we have
$$\Psi[x_{\a+\sg},x_{-\a}]=\eta_{\a+\sg}\eta_{-\a}[x_{\a+\sg},x_{-\a}]=\eta_\sg[x_{\a+\sg},x_{-\a}].$$ Since $E_c\cap E_\sg$ is spanned by $1$-dimensional spaces $[E_{\a+\sg}, E_\a]$, $\a\in R^\times$, we get
\begin{equation}\label{diag3}\Psi_{|_{E_c\cap E_\a}}=\eta_\sg\id_{|_{E_c\cap E_\a}},\qquad(\a\in R).\end{equation}

Now, consider an arbitrary element $x_\sg\in E_\sg$, $\sg\in R^0$. Then for each $\a\in R^\times$ such that $\a+\sg\in R$ we have
$$\eta_{\a+\sg}[x_\a,x_\sg]=\Psi[x_\a,x_\sg]=\eta_{\a}[x_\a,\Psi(x_\sg)].$$
So $[x_\a,\eta_\sg x_\sg-\Psi(x_\sg)]=0$. We note that this equality holds even if  $\a+\sg\not\in R$, since in this case $[x_\a,x_\sg]=0=[x_\a,\Psi(x_\sg)].$ Thus $\eta_\sg x_\sg-\Psi(x_\sg)$ is in the centralizer of $E_c$. Since $E$ is tame, $\eta_\sg x_\sg-\Psi(x_\sg)\in E_c\cap E_\sg$.
Then by (\ref{diag3}), $\Psi(\eta_\sg x_\sg-\Psi(x_\sg))=\eta_\sg(\eta_\sg x_\sg-\Psi(x_\sg)).$
Therefore, $\eta_\sg^{-1}\Psi_{|_{E_\sg}}$ satisfies the equation $(t-1)^2=0$. Since $\Psi$ is of finite order, we get $\Psi(x_\sg)=\eta_\sg x_\sg$ for each $x_\sg\in E_\sg$, and so
$\Psi_{|_{E_\a}}=\eta_\a\id_{|_{E_\a}},$ $\a\in R$, as required.
\qed

\begin{DEF}\label{aug31a}
Let $\Psi$ be a finite order diagonal automorphism for $E$. For $\sg\in R^0$, we set $\eta_\sg:=\eta_{\a+\sg} \eta_{-\a}$, for some $\a\in R^\times.$ By Lemma \ref{bims30}, $\eta_\sg$ is well-defined.
Proposition \ref{aug30c} below shows that the map $\psi:R\rightarrow\bbbk^\times$ defined by $\psi(\a)=\eta_\a$, $\a\in R$, is a core-character, hereafter, we refer to it as the {\it core-character associated to} $\Psi$. 
\end{DEF}

\begin{pro}\label{aug30c}
	Let $\Psi$ be a finite-order diagonal automorphism. The map
	 $\psi:R\rightarrow\bbbk^\times$ defined by $\psi(\a)=\eta_\a$, $\a\in R$ is a finite-order core-character for $R$, satisfying 
	 $\psi(-\a)=\psi(\a)^{-1}$, $\a\in R$.
\end{pro}

\proof Let $\psi$ be the map associated to $\Psi$, given in  Definition \ref{aug31a}. Since $\Psi$ stabilizes pointwise $\hh$, we note that if  $0\not=x_{\pm\a}\in E_{\pm\a}$, $\a\in R^\times$, we have
$$[x_\a,x_{-\a}]=\Psi[x_\a,x_{-\a}]=\eta_\a\eta_{-\a}[x_\a,x_{-\a}],$$
implying that $\psi(-\a)=\psi(\a)^{-1}$. Also if $\sg\in R^0$ and $\a,\a+\sg\in R^\times$, then
$$\psi(\sg)\psi(-\sg)=\eta_\sg\eta_{-\sg}=\eta_{\a+\sg}\eta_{-\a}\eta_{-\a-\sg}\eta_{\a}=1.
$$

Next, we show that $\psi$ is a core-character. If $\a,\b,\a+\b\in R^\times$ and $x_\a,x_\b$ are non-zero elements of $E_\a$ and $E_\b$, then $[x_\a,x_\b]\not=0$, and
$$
\eta_{\a+\b}[x_\a,x_\b]=\Psi[x_\a,x_\b]=\eta_\a\eta_\b[x_\a,x_\b].
$$
Thus $\psi(\a+\b)=\eta_{\a+\b}=\eta_\a\eta_\b=\psi(\a)\psi(\b).$
Also if $\b\in R^0,\a,\a+\b\in R^\times$, then by the way  $\eta_\b$ is defined, we have
$$\psi(\a+\b)=\eta_{\a+\b}=\eta_{\a+\b}\eta_{-\a}\eta_{\a}=\eta_\b\eta_\a=\psi(\a)\psi(\b).$$ Thus $\psi$ is a core-character. Since $\Psi$ is of finite-order, so is $\psi$. 
\qed

We proceed with stating a few lemmas which will be used frequently in the sequel.
\begin{lem}\label{bim4} Let $\sg\in R^0$.
	
	(i) If $\a,\b,\a+\sg,\b+\sg,\a-\b\in R^\times$, then
	$\eta^\a_\sg=\eta^\b_\sg$.
	
	(ii)  If $\a,\b,\a+\b,\a+\b+\sg,\b+\sg\in R$, then
	$\eta^{\a+\b}_\sg=\eta^\b_\sg$.
\end{lem}

\proof (i) Applying 
$\Psi$ to the non-zero element 
$[x_{\a+\sg},x_{-\b-\sg}]$, gives $\eta_{\a+\sg}\eta_{-\b-\sg}=\eta_{\a-\b}=\eta_{\a}\eta_{-\b}$, as required.

(ii) Set $\gamma=\a+\b$. Then the elements $\gamma$, $\b$ and $\sg$ satisfy conditions of part (i),  with $\gamma$ in place of $\a$.\qed


\begin{lem}\label{Aug23b}
	Let $\sg\in R^0$ and $\a,\b,\a+\sg,\b+\sg\in R^\times$. For
	$\gamma\in\{\pm\a,\pm\b,\pm(\a+\sg)\}$, fix $0\not=x_\gamma\in E_\gamma$. If $\b+\a+\sg\not\in R$ but $\b-\a\in R^\times$, then
	$[x_\b,[x_{\a+\sg},x_{-\a}]]\not=0$.
\end{lem}

\proof By assumptions, $[x_{\a+\sg},[x_{\b},x_{-\a}]]\not=0$ and
$[x_{-\a},[x_{\a+\sg},x_{\b}]]=0$. Therefore, we are done by the Jacobi identity.\qed	


\begin{lem}\label{Aug23a}
Let  $\sg\in R^0$ and {$\a,\a+\sg\in R^\times$}. Then

(i) 	$[E_\a,[E_{\a},E_{\sg-\a}]]\not=\{0\}$.

(ii) $\eta_\a^2=\eta_{\a+\sg}\eta_{\a-\sg}$.
\end{lem}

\proof (i)  Consider the $\a$-string through $\sg$ (see (\ref{string})), namely
$$\{\sg-d\a,\ldots,\sg,\ldots\sg+u\a\}=R\cap\{\sg+n\a\mid n\in\bbbz\},\quad(u,d\geq 0).$$ Since $(\a,\sg)=0$, $1\leq u=d\leq 2$.
Fix $e_{\pm\a}\in E_{\pm\a}$ with $[e_\a,e_{-\a}]=h_\a$, and let $0\not= x_{\sg-d\a}\in E_{\sg-d\a}.$ Then $[e_{-\a},x_{\sg-d\a}]=0$. Therefore by \cite[Lemma I.1.21]{AABGP97},
$$\ad(e_\a)^N(x_{\sg-d\a})\not=0\andd \ad(e_\a)^{N+1}(x_{\sg-d\a})=0,$$
where $N:=(\sg-d\a,-\a^\vee)=2d\geq 2$.
Then, $z:=[e_{\a},[e_{\a},(\ad e_\a)^{d-1}(x_{\sg-d\a})]]\not=0$, as
$d+1\leq 2d=N$. Now $z\in[E_\a,[E_\a,E_{\sg-\a}]]$ and we are done.

(ii)  Applying $\Psi$ to the element $[x_\a,[x_\a,x_{-\a+\sg}]$ which is non-zero by part (i), we get $\eta_\a\eta_\a\eta_{-\a+\sg}=\eta_{\a+\sg}$, as required.\qed

\subsection{Simply laced Cartan automorphisms}\label{bims22}
We study the class of finite-order Cartan automorphisms, depending on type $X$ of $E$. We begin with simply laced types.

\begin{pro}\label{lemaug31a} Let $X$ be simply laced of rank $>1$, and $\Psi$ a finite order Cartan automorphism for $E$. Then

	(i) $\Psi$ is diagonal,
	
	(ii) the core-character $\psi$ associated to $\Psi$ is a character.
\end{pro}

\proof Under the assumptions, we know from \pref{bims23} that $S=\Lam=R^0$ is a lattice, and
$R^\times=\rd^\times+\Lam$.
By Lemma  \ref{bims30}, it is enough to show that for $
\lam\in R^0$ and $\a,\b\in R^\times$, $\eta^\a_\lam=\eta^\b_\lam$. We first claim that, 
\begin{equation}\label{tempaug31-1}
\eta^{\dot\a}_\lam=	\eta^{\dot\b}_\lam,\quad
\dot\a,\dot\b\in \dot R^\times, \lam\in R^0.	
\end{equation}
To prove this, we primarily note that if $\dot\a,\dot\b,
\dot\a+\dot\b\in\rd^\times$, then by
applying $\Psi$ on the non-zero element
$[x_{-\dot\a},[x_{\dot\a+\lam}, x_{\dot\b}]]$, we obtain
$\eta_{\dot\b+\lam}=\eta_{-\dot\a}\eta_{\dot\a+\lam}\eta_{\dot\b}$, or equivalently
\begin{equation}\label{fer1}
\eta^{\dot\a}_\lam=\eta^{\dot\b}_\lam,\quad(\dot\a,\dot\b,\dot\a+\dot\b\in\rd^\times).
\end{equation}
Also, we see from Lemma \ref{Aug23a}(ii) that
\begin{equation}\label{bim1}
\eta_\lam^{\dot\a}=
\eta_\lam^{-\dot\a},\quad(\dot\a\in\rd^\times,\lam\in R^0).
\end{equation}

Now let $\dot\a,\dot\b\in \rd^\times$. By (\ref{fer1}), we may assume that $\dot\a\not=\pm\dot\b$. Since $R$ is irreducible, we can find a sequence $\dot\a_1=\dot\a,\dot\a_2,\ldots\,\dot\a_{n-1},\dot\a_n=\dot\b$ of roots in $\rd^\times$
with $(\dot\a_i,\dot\a_{i+1})\not=0$, for all $i$. We may assume also that 
$\dot\a_i\not=\pm\dot\a_{i+1}.$ By the root string property, either $\dot\a_i+\dot\a_{i+1}\in\rd^\times$ or 
${\dot\a}_i-{\dot\a}_{i+1}\in\rd^\times$. Therefore, by changing signs if necessary, we may assume that ${\dot\a}_i+{\dot\a}_{i+1}\in\rd^\times$
for all $i$. Then $\eta_\lam^{\dot\a}{=\eta_\lam^{\dot\a_1}=}\eta_{\lam}^{\dot\a_2}=\cdots=\eta_{\lam}^{\dot\a_n}=\eta_\lam^{\dot\b}.$
This proves (\ref{tempaug31-1}). 

Next let $\a\in R^\times$, $\lam\in R^0$. Then
$\a=\dot\a+\sg$ for some $\dot\a\in\rd$, $\sg\in R^0$. Since $\rank R >1$, we can choose $\dot\b\in\rd^\times$ such that $\dot\a+\dot\b\in R^\times$.
Then
\begin{eqnarray*}
	\eta^\a_\lam
	&=&
	\eta_{\dot\a+\sg+\lam}\eta_{-\dot\a-\sg}\\
	&=&\eta_{(\dot\a+\dot\b+\sg)+(-\dot\b+\lam)}\eta_{-\dot\a-\sg}\\
	&=&\eta_{\dot\a+\dot\b+\sg}\eta_{-\dot\b+\lam}\eta_{-\dot\a-\sg}\\
	&=&\eta_{\dot\a+\sg}\eta_{\dot\b}\eta_{-\dot\b+\lam}\eta_{-\dot\a-\sg}\\
	&=&\eta^{-\dot\b}_\lam\\
	(\hbox{by (\ref{tempaug31-1}}) &=&\eta^{\dot\a}_\lam.
\end{eqnarray*}
Therefore, we see from this by again applying (\ref{tempaug31-1}) that  if
$\b=\dot\b+\sg'\in R^\times$, $\dot\b\in\rd^\times$, $\sg'\in R^0$, we have
$\eta^\a_\lam=\eta^{\dot\a}_\lam=\eta^{\dot\b}_\lam=\eta^\b_\lam.$
This completes the proof of (i).

(ii) By part (i) and Proposition \ref{aug31a}, $\psi$ is a core-character.
To see that $\psi$ is a character, let $\sg,\sg'\in S=R^0$. Then for $\dot\a\in\rd^\times$ we have
\begin{eqnarray*}
	\psi(\sg+\sg')
	&=&\psi(\dot\a+\sg-\dot\a+\sg')\\
	&=&
	\psi(\dot\a+\sg)\psi(-\dot\a+\sg')\\
	&=&\psi(\dot\a)\psi(\sg)\psi(-\dot\a)\psi(\sg')\\
	&=&\psi(\sg)\psi(\sg').
\end{eqnarray*}
\qed

\subsection{Non-simply laced Cartan automorphisms}\label{bims22}
We continue with the study of Cartan automorphisms, for non-simply laced extended affine Lie algebras. We assume that $E$ is a non-simply laced reduced extended affine Lie algebra of type $X$, with root system $R$. We recall from \pref{bims23} that 
$R=R^0\cup(\rds+S)\cup(\rdl+L)$, with $R^0=(S+S)$. Also, we recall that
$\rd=\rds\cup\rdl\cup\{0\}$ is an irreducible finite root system which has one of the types $X=B_\ell,$ $C_\ell$, $F_4$ or $G_2$. We denote the highest short and long roots of $\rd$ with $\dot\theta_s$ and $\dot\theta_\ell$, respectively.  As before $\Psi$ is a finite order Cartan automorphism.

\begin{lem}\label{Sept1-a} We have
	\begin{equation}\label{test3}
	\eta^{\dot\theta_\ell}_\lam=	\eta^{\dot\theta_s}_\lam,\quad (\lam\in L),
\end{equation}
	and
\begin{equation}\label{bim5}
	\eta^{\dot\a}_\sg=\eta^{\dot\a'}_\sg,\quad(\dot\a,\dot\a'\in\rds,\sg\in S).
\end{equation}

	\end{lem}

\proof 
Let $\lam\in L$. Take $\a=\dot\theta_\ell+\lam$ and $\b=\dot\theta_s+\lam$. Then $\a-\b=\dot\theta_\ell-\dot\theta_s$, which is a root.Thus  applying Lemma \ref{bim4} we get (\ref{test3}).

Next, let $\dot\a,\dot\a'\in\rds$ and $\sg\in S$. If $X$ is of type $B_\ell$, then
$\dot\a\pm\dot\a'\in\rdl$ and so we are done by Lemma \ref{bim4}. If $X=C_\ell$ $(\ell\geq 3)$, $F_4$ or $G_2$, then by \cite[Theorem II.2.27]{AABGP97}, $R':=S\cup(\rds+S)$ is an extended affine root system of simply laced type, so the Proposition \ref{lemaug31a} applies. 
\qed


\begin{lem}\label{lemsept1a}
 (i) If $X=B_\ell$, then
	$\eta^\a_\lam=\eta^{\b}_\lam$, $\a,\b\in R^\times$, $\lam\in L$.
	
	(ii)	If  $X=C_\ell,\ell\geq 3$, $F_4$ or $G_2$, then
	$\eta^\a_\sg=\eta^{\b}_\sg$, $\a,\b\in R^\times$, $\sg\in R^0$.
	\end{lem}

\proof (i) 
First, we claim that
\begin{equation}\label{bim9}
	\eta^\a_\lam=\eta^\b_\lam,\quad(\lam\in L,\a,\b,\a+\lam,\b+\lam\in R_{\lg}).
\end{equation}
To see this, we consider two cases $\ell\geq 3$ and $\ell=2$, separately.
Assume first that $\ell\geq 3$. Then $L$ is a lattice and so $\a+\lam\in R^\times$ for any $\a\in R^\times$ and $\lam\in L$.
We also note that by  \cite[Theorem II.2.37]{AABGP97}, $R':=L\cup(\rdl+L)$ is
a simply laced extended affine root system of $\rank >1$. 
Therefore, by Proposition \ref{lemaug31a}(i), the equality (\ref{bim9}) holds.

Next, assume that $\ell=2$. Let $\lam\in L$, $\b\in R_{lg}$ and $\b+\lam\in R$. Then $\b=\dot\b+\lam_\ell$, for some $\dot\b\in\rdl$ and some $\lam_\ell\in L$. From the realization of $\rd$, we know that $\dot\b=\dot\a_1+\dot\a_2$ for some short roots $\dot\a_1,\dot\a_2\in\rds$. Since $\lam\in L\sub S$, we have
$[x_{\dot{\a_1}+\dot{\a_2}+\lam_\ell+\lam},
[x_{-\dot{\a_1}-\dot{\a_2}},x_{\dot\a_1-\lam_\ell}]]\not=0$. 
Applying $\Psi$ to this, we get
$$\eta_{\dot\a_1+\dot\a_2+\lam_\ell+\lam}\eta_{-\dot\a_1-\dot\a_2}\eta_{\dot\a_1-\lam_\ell}=\eta_{\dot\a_1+\lam},$$ or equivalently,
 $$\eta_{\dot\a_1+\dot\a_2+\lam_\ell+\lam}\eta_{-\dot\a_2-\lam_\ell}=\eta_{\dot\a_1+\lam}.$$
 Now inserting $\eta_{-\dot\a_1}$ to the both sides of this last equality we get
 $\eta^\b_\lam=\eta^{\dot\a_1+\dot\a_2+\lam_\ell}_\lam=\eta^{\dot\a_1}_\lam$, and so the claim (\ref{bim9}) holds by (\ref{bim5}).  
 
 Now combining (\ref{bim9}) with Lemma \ref{Sept1-a} gives
\begin{equation}\label{bim7}
\eta^{\gamma}_\lam=\eta_\lam^{\dot\a},\quad(\dot\a\in\rd^\times, \gamma\in R_{lg},\lam\in L).
\end{equation}

Next, let $\a\in R_{sh}$ and $\lam\in L$. Then $\a=\dot\a_s+\lam_s$ for some $\dot\a_s\in\rds$ and $\lam_s\in S$. Choose $\dot\b_s\in\rds$ such that $\dot\a_s-\dot\b_s\in\rdl$. Then
\begin{eqnarray*}
\eta_{\lam}^{\a}
	&=&
	\eta_{\dot\a_s+\lam_s+\lam}\eta_{-\dot\a_s-\lam_s}\\
	&=&
	\eta_{\dot\a_s+\dot\b_s-\dot\b_s+\lam_s+\lam}\eta_{-\dot\a_s-\lam_s}\\
	&=&
	\eta_{(\dot\a_s-\dot\b_s)+\lam}\eta_{\dot\b_s+\lam_s}\eta_{-\dot\a_s-\lam_s}\\
	&=&
	\eta_{(\dot\a_s-\dot\b_s)+\lam}\eta_{-\dot\a_s+\dot\b_s}\eta_{\dot\a_s-\dot\b_s}\eta_{\dot\b_s+\lam_s}\eta_{-\dot\a_s-\lam_s}\\
	&=&
	\eta^{\dot\a_s-\dot\b_s}_\lam\eta_{-\dot\b_s}\eta_{\dot\b_s+\lam_s}\eta_{-\dot\a_s-\lam_s}\eta_{\dot\a_s}\\
	&=&
	\eta^{\dot\a_s-\dot\b_s}_\lam\eta^{-\dot\a_s}_{-\lam_s}\eta^{\dot\b_s}_{\lam_s}.
\end{eqnarray*}
{Now from (\ref{bim5}), we have $\eta_{-\lam_s}^{-\dot\a_s}=\eta^{\dot\b_s}_{-\lam_s}$. Therefore,
$\eta_{\lam}^\a=	\eta^{\dot\a_s-\dot\b_s}_\lam\eta^{\dot\b_s}_{-\lam_s}\eta^{\dot\b_s}_{\lam_s}=	
	\eta^{\dot\a_s-\dot\b_s}_\lam.$
This together with (\ref{bim7}) completes the proof of (i).}

(ii) First we consider type $X=C_\ell$, $\ell\geq 3$. In this case $R^0=S$ is a lattice.
Now by \cite[Theorem II.2.37]{AABGP97}, $R'=R^0\cup(\rds+S)$ is a simply laced extended affine root system. In fact if $\ell=3$, then $R'$ is of type $A_3$ and otherwise is of type $D_\ell$. Then by 
Proposition \ref{lemaug31a}(i), 
\begin{equation}\label{temp1}
	\eta_\sg^\a=\eta_\sg^\b,\qquad(\sg\in R^0,
\a,\b\in R_{sh}=\rds+S).
\end{equation}

Next, let $\sg\in R^0=S$, $\gamma,\gamma+\sg\in R_{\lg}=\rdl+L$.
Then $\gamma=\dot\gamma+\lam$ for some $\dot\gamma\in\rdl$, $\lam\in L$. 
Since  $\gamma+\sg\in R$, we have $\lam+\sg\in L$. From realizations of finite root systems of type $C_\ell$, we know that there exist roots,
$\dot\a_1,\dot\a_2\in\rds$ such that $\dot\a_1+\dot\a_2=\dot\gamma$.
Therefore,
\begin{eqnarray*}
	\eta_\sg^{\gamma}
	&=&
	\eta_{\dot\gamma+\lam+\sg}\eta_{-\dot\gamma-\lam}\\
	&=&
	\eta_{(\dot\a_1+\sg)+(\dot\a_2+\lam)}\eta_{-\dot\a_1+(-\dot\a_2-\lam)}\\
	&=&
	\eta_{\dot\a_1+\sg}\eta_{\dot\a_2+\lam}\eta_{-\dot\a_1}\eta_{-\dot\a_2-\lam}\\
	&=&
	\eta^{\dot\a_1}_\sg.
	\end{eqnarray*}
This together with (\ref{temp1}) completes the proof for type $C_\ell$.

Finally, let $X=F_4$ or $G_2$, and $\sg\in R^0$. In this case both $S$ and $L$ are lattices.
Moreover, both $R'=S\cup(\rds+S)$ and $R''=L\cup(\rdl+L)$ are simply laced extended affine root systems of rank $>1$. Therefore, by Proposition \ref{lemaug31a},
\begin{equation}\label{temp2}
\eta_\sg^\a=\eta_\sg^\b,\quad(\a,\b\in R_{sh}, \sg\in S),\hbox{ or }(\a,\b\in R_{\lg},\sg\in L).
\end{equation}
Now if $\sg\in S=R^0$, $\a\in R_{sh}, \b\in R_{lg}$ and $\a+\sg,\b+\sg\in R$, then since $L$ is a lattice, we get $\sg\in L$, so
by (\ref{temp2}) and (\ref{test3}),
$$\eta_\sg^\b=\eta_\sg^{\dot\theta_l}=\eta_\sg^{\dot\theta_s}=\eta_\sg^{\a}.$$
This completes the proof.\qed

\begin{cor}\label{red1}
For $\sg\in S$ and $\dot\a,\dot\a'\in\rd^\times$, $\eta^{\dot\a}_\sg=\eta^{\dot\a'}_\sg$.
\end{cor}

\proof Combining  (\ref{bim5}) with Lemma \ref{lemsept1a}, we get the result.\qed

\begin{pro}\label{bim11}
	Any finite-order Cartan automorphism $\Psi$ of type $C_\ell$, $\ell\geq 3$, $F_4$ or $G_2$, is diagonal, and the corresponding map $\psi:R\rightarrow \bbbk^\times$ is a character.
\end{pro}
\proof By Lemma \ref{lemsept1a}, $\Psi$ is diagonal. From Proposition \ref{aug30c}, we know that $\psi$ is a core-character. Now any two isotropic roots $\sg_1,\sg_2\in R^0=S$ can be considered as roots of the simply laced extended affine root system  $R':=S\cup(\rds+S)$.  Since $\rank R'>1$, we have from Proposition \ref{lemaug31a} that $\psi(\sg_1+\sg_2)=\psi(\sg_1)\psi(\sg_2)$. Thus $\psi$ is a character.
\qed

\subsection{Cartan automorphisms of Type $B_\ell $, $\ind(R)=0$}\label{type b} 	
Throughout this section $R$ is an extended affine root system of type $B_\ell$, with root system $R$, such that $\ind(R)=0$.
Let $t$ be the twist number of $R$.
From \pref{july20b}, we see that 
$S=S_1\oplus \Lam_2$, $L=2 \Lam_1\oplus S_2$, where $\Lam_1=\la S_1\ra$ and $\Lam_2=\la S_2\ra.$ Moreover, we have $\ind(S_1)=t$, and $S_1$ can be written in the form $S_1=\uplus_{i=1}^t(\tau_i+2\Lam_1)$, where
$\Lam_1=\la S_1\ra=\sum_{i=1}^t\bbbz\tau_i$, and $\tau_i$'s are distinct modulo $2\Lam_1$. For the rest of this section, we fix $\tau_i$'s as above. Recall that $\Psi$ is a finite order Cartan automorphism. 
	
\begin{lem}\label{b100}
	Let $\a,\in R_{sh}$ and $\lam_1,\ldots,\lam_k\in L$, then
	$$\eta^\a_{\lam_1+\cdots+\lam_k}=\eta^\a_{\lam_1}\cdots\eta^\a_{\lam_k}
	=\eta^{\dot\theta_s}_{\lam_1}\cdots\eta^{\dot\theta_s}_{\lam_k}=\eta^{{\dot\theta}_s}_{\lam_1+\cdots+\lam_k}.$$	In particular, 	$\eta^\a_{\lam+\lam'}=\eta^\a_{\lam}\eta ^\a_{\lam'}=\eta^{\a'}_{\lam}\eta ^{\a'}_{\lam'}=\eta^{\a'}_{\lam+\lam'}$, $\a,\a'\in R_{sh}$, $\lam,\lam'\in\la L\ra$.
\end{lem}

\proof  The proof will make frequent use of the fact that $R_{sh}+\la L
\ra\sub R_{sh}$, without explicitly stating it every time.
We have $\a=\dot\a+\sg$, where $\dot\a\in\rds$ and $\sg\in S$. Take $\dot\b\in\rds$ with $\dot\b\not=\pm\dot\a$. Then $\dot\a\pm\dot\b\in\rdl$.
Now
\begin{eqnarray*}
	\eta^\a_{\lam_1+\lam_2}&=&\eta_{\dot\a+\sg+\lam_1+\lam_2}\eta_{-\dot\a-\sg}\\
	&=&\eta_{\dot\a+\dot\b-\dot\b+\sg+\lam_1+\lam_2}\eta_{-\dot\a-\sg}\\
	&=&\eta_{\dot\a-\dot\b+\lam_1}\eta_{\dot\b+\sg+\lam_2}\eta_{-\dot\a-\sg}\\
	&=&\eta_{\dot\a-\dot\b+\lam_1}\eta_{\dot\b-\dot\a+\lam_2}\\
	&=&\eta_{\dot\a+\lam_1}\eta_{-\dot\b}\eta_{\dot\b}\eta_{-\dot\a+\lam_2}\\
		&=&\eta_{\dot\a+\lam_1}\eta_{-\dot\a}\eta_{\dot\a}\eta_{-\dot\a+\lam_2}\\
	&=&\eta^{\dot\a}_{\lam_1}\eta^{-\dot\a}_{\lam_2}\\
	(\hbox{by Lemma \ref{lemsept1a}(i)})&=&\eta^{{\dot\theta}_s}_{\lam_1}\eta^{{\dot\theta}_s}_{\lam_2}.
\end{eqnarray*}
Now using induction on $k$, we have
\begin{eqnarray*}
	\eta^\a_{\lam_1+\cdots+\lam_k}&=&\eta_{\dot\a+\sg+\lam_1+\cdots+\lam_k}\eta_{-\dot\a-\sg}\\
	&=&\eta_{\dot\a-\dot\b+\lam_1}\eta_{\dot\b+\sg+\lam_2+\cdots+\lam_k}\eta_{-\dot\a-\sg}\\
	&=&\eta_{\dot\a-\dot\b+\lam_1}\eta_{\dot\b+\sg+\lam_2+\cdots+\lam_k}\eta_{-\dot\b-\sg}\eta_{\dot\b+\sg}\eta_{-\dot\a-\sg}\\
	&=&\eta^{\dot\a-\dot\b}_{\lam_1}\eta_{\dot\a-\dot\b}\eta^{\dot\b+\sg}_{\lam_2+\cdots+\lam_k}\eta_{\dot\b+\sg}\eta_{-\dot\a-\sg}\\
	(\hbox{induction steps})&=&\eta^{\dot\theta_s}_{\lam_1}\eta^{\dot\theta_s}_{\lam_2}\cdots\eta^{\dot\theta_s}_{\lam_k}\eta^{\dot\b}_{\sg}\eta^{-\dot\a}_{-\sg}\\
	(\hbox{by Corollary \ref{bim5}})	&=&\eta^{\dot\theta_s}_{\lam_1}\eta^{\dot\theta_s}_{\lam_2}\cdots\eta^{\dot\theta_s}_{\lam_k}\eta^{\dot\a}_{\sg}\eta^{\dot\a}_{-\sg}\\
	&=&\eta^{\dot\theta_s}_{\lam_1}\eta^{\dot\theta_s}_{\lam_2}\cdots\eta^{\dot\theta_s}_{\lam_k}.
\end{eqnarray*}
This completes the proof.\qed

\begin{cor}\label{aban1}
	Let $\a,\b\in R_{sh}$, $\lam\in L$ and $\sg\in R^0$. Then
	$$\eta^\a_{\lam+2\sg}=\eta^\b_{\lam+2\sg}\Longleftrightarrow \eta^\a_\lam=\eta^\b_\lam.
	$$
	\end{cor}

\proof First, we note that  $\a+\lam+2\sg\sub R_{sh}+L+2\Lam\sub R_{sh}$. By Lemma \ref{b100}, we have
$\eta^\a_{\lam+2\sg}=\eta^\a_\lam\eta^\a_{2\sg}$, and $\eta^\b_{\lam+2\sg}=\eta^\b_\lam\eta^\b_{2\sg}$. So, we are done by Lemma \ref{lemsept1a}(i),
$\eta^\a_{2\sg}=\eta^\b_{2\sg}$.\qed

\begin{lem}\label{red2}
	For $\sg\in S$ and $\a,\b\in R^\times$, $\eta^\a_\sg=\eta^\b_\sg$.
\end{lem}
\proof
We begin by proving that
\begin{equation}\label{red3}
\eta^\a_\sg=\eta^\b_\sg,\quad (\sg\in S,\a,\b\in R_{sh}).
\end{equation}
First, assume that $\sg\in S_1$. Let $\a\in R_{sh}$. Using Corollary \ref{aban1}, we may assume that $\a=\dot\a+\ep_i\tau_i+\lam_2$, and
$\sg=\tau_j$, for some $\dot\a\in\rds$, $\lam_2\in\Lam_2$, $1\leq i,j\leq t$, $\ep_i\in\{0, 1\}$. If $\a+\sg\in R$, then $\ep_i\tau_i+\tau_j+\lam_2\in S$. Since $\ind(R)=0$, 
either $i=0$, or $\ep_i=1$ and  $i=j$. In the latter case, we have
\begin{eqnarray*}
	\eta^\a_\sg&=&\eta_{\dot\a+2\tau_i+\lam_2}\eta_{-\dot\a-\tau_i-\lam_2}\\
	&=&\eta_{\dot\a+2\tau_j+\lam_2}\eta_{-\dot\a-2\tau_j}\eta_{\dot\a+2\tau_j}\eta_{-\dot\a-\tau_j}\eta_{\dot\a+\tau_j}\eta_{-\dot\a-\tau_j-\lam_2}\\
	&=&\eta^{\dot\a+2\tau_j}_{\lam_2}\eta_{\dot\a+2\tau_j}\eta^{-\dot\a-\tau_j}_{\lam_2}\eta_{-\dot\a-\tau_j}\\
(\hbox{by Lemma \ref{lemsept1a}(i)})	&=&\eta^{\dot\a}_{\lam_2}\eta_{\dot\a+2\tau_j}\eta_{-\dot\a-\tau_j}\eta^{-\dot\a}_{-\lam_2}\\
	&=&\eta_{\dot\a+2\tau_j}\eta_{-\dot\a-\tau_j}\\
	(\hbox{by Lemma \ref{Aug23a}})&=&
	(\eta_{\dot\a+\tau_j})^2\eta_{-\dot\a}\eta_{-\dot\a-\tau_j}\\
	&=&\eta^{\dot\a}_{\tau_j}.
\end{eqnarray*}
If $\ep_i=0$, then
\begin{eqnarray*}
	\eta^\a_\sg&=&\eta_{\dot\a+\tau_j+\lam_2}\eta_{-\dot\a-\lam_2}\\
	&=&\eta_{\dot\a+\tau_j+\lam_2}\eta_{-\dot\a-\tau_j}\eta_{\dot\a+\tau_j}\eta_{-\dot\a-\lam_2}\\
	&=&\eta^{\dot\a+\tau_j}_{\lam_2}\eta_{\dot\a+\tau_j}\eta_{-\dot\a}\eta_{\dot\a}\eta_{-\dot\a-\lam_2}\\
(\hbox{by Lemma \ref{lemsept1a}(i)})	&=&\eta^{\dot\a}_{\lam_2}\eta^{\dot\a}_{\tau_j}\eta^{-\dot\a}_{-\lam_2}\\
	&=&\eta^{\dot\a}_{\tau_j}.
\end{eqnarray*}
Thus in either cases  $\eta^\a_{\sg}=\eta^{\dot\a}_{\sg}$, and so by Corollary \ref{red1}, 
$\eta^\a_\sg=\eta^\b_\sg$, for any
$\a,\b\in R_{sh}$, $\sg\in S_1$.

Generally,  assume that $\sg$ is an arbitrary element of $S$. Then $\sg=\sg_1+\lam_2$, where $\sg_1\in S_1$ and $\lam_2\in\Lam_2$.
Then
\begin{eqnarray*}
	\eta^\a_\sg&=&\eta_{\a+\sg_1+\lam_2}\eta_{-\a}\\
	&=&\eta_{\a+\sg_1+\lam_2}\eta_{-\a-\sg_1}\eta_{\a+\sg_1}\eta_{-\a}\\
	&=&\eta^{\a+\sg_1}_{\lam_2}\eta^{\a}_{\sg_1}\\
	&=&\eta^{\a}_{\lam_2}\eta^{\a}_{\sg_1}.
\end{eqnarray*}
Thus by Lemma \ref{lemsept1a}(i),
\begin{equation}\label{red4}
\eta^\a_{\sg_1+\lam_2}=\eta^\a_{\sg_1}\eta^\a_{\lam_2}=\eta^{\b}_{\sg_1}\eta^\b_{\lam_2}=\eta^\b_{\sg_1+\lam_2},\quad(\a,\b\in R_{sh},\sg_1\in S_1,\lam_2\in\Lam_2)
\end{equation}
This completes the proof of (\ref{red3}).

Next, let $\b\in R_{lg}$ and $\sg\in S$ with $\b+\sg\in R$. By Corollary \ref{aban1}, we may assume that $\b=\dot\b+\sg_2$, and
$\sg=\sg_1+\lam_2$, where $\dot\b\in\rdl$, $\sg_1\in S_1$, $\sg_2\in S_2$ and $\lam_2\in\Lam_2$. Since $\a+\sg\in R$, we
get $\sg_1+\sg_2+\lam_2\in L\sub 2\Lam_1\oplus\Lam_2$, and so $\sg_1\in 2\Lam_1$. Also $\sg_2\in L+L+2\Lam_2\sub\la L\ra$.  Take $\dot\a,\dot\gamma\in\rds$ such that $\dot\b=\dot\a+\dot\gamma$. Then
\begin{eqnarray*}
	\eta^\b_\sg&=&\eta_{\dot\a+\dot\gamma+\sg_1+\sg_2+\lam_2}\eta_{-\dot\a-\dot\gamma-\sg_2}\\
	&=&\eta_{\dot\a+\sg_1+\lam_2}\eta_{\dot\gamma+\sg_2}\eta_{-\dot\a}\eta_{-\dot\gamma-\sg_2}\\
	&=&\eta_{\dot\a+\sg_1+\lam_2}\eta_{-\dot\a-\sg_1}\eta_{\dot\a+\sg_1}\eta_{-\dot\a}\\
	&=&\eta^{\dot\a+\sg_1}_{\lam_2}\eta^{\dot\a}_{\sg_1}\\
(\hbox{by (\ref{red4})})	&=&\eta^{\dot\a}_{\lam_2}\eta^{\dot\a}_{\sg_1}=\eta^{\dot\a}_{\sg}.
\end{eqnarray*}
Since $\dot\a$ is a short root, the proof is completed by (\ref{red3}).\qed

\begin{DEF}\label{bim10}
  For $\sg\in S$, we define $\eta_\sg=\eta^\a_\sg$, for some $\a\in R^\times.$ By Lemma \ref{red2}, this is well-defined.
	\end{DEF}


\begin{lem}\label{bim15} (i) $\eta_{\lam+2\sg}=\eta_\lam\eta_{2\sg}$, $\lam\in S$, $\sg\in R^0$.
	
	(ii) $\eta_{\lam_1+\lam_2}=\eta_{\lam_1}\eta_{\lam_2}$, $\lam_1\in S_1$, $\lam_2\in\la S_2\ra$.
\end{lem}

\proof Take $\dot\a_1,\dot\a_2\in\rds$ with $\dot\a_1\not=\pm\dot\a_2$.
Then $\dot\a_1\pm\dot\a_2\in\rdl$.

(i) For $\sg\in R^0$, we have
$2\sg\in 2R^0\sub 2\la S\ra\sub L\sub S$. Now
\begin{eqnarray*}
	\eta_{\lam+2\sg}&=&\eta_{\dot\a_1+\lam+2\sg}\eta_{-\dot\a_1}\\
	&=&\eta_{(\dot\a_2+\lam)+(\dot\a_1-\dot\a_2+2\sg)}\eta_{-\dot\a_1}\\
	&=&\eta_{\dot\a_2+\lam}\eta_{\dot\a_1-\dot\a_2+2\sg}\eta_{-\dot\a_1}\\
	&=&(\eta_{\dot\a_2+\lam}\eta_{-\dot\a_2})\eta_{\dot\a_2}\eta_{\dot\a_1-\dot\a_2+2\sg}\eta_{-\dot\a_1}\\
	&=&\eta_\lam\eta_{\dot\a_1-\dot\a_2+2\sg}\eta_{\dot\a_2-\dot\a_1}\\
	&=&\eta_\lam\eta_{2\sg}.
\end{eqnarray*}

(ii) We have
\begin{eqnarray*}	\eta_{\lam_1+\lam_2}
	&=&
	\eta^{\dot\a_1}_{\lam_1+\lam_2}\\
	(\hbox{by Lemma \ref{b100}})
	&=& \eta^{\dot\a_1}_{\lam_1}\eta^{\dot\a_1}_{\lam_2}\\
	&=&\eta_{\lam_1}\eta_{\lam_2}.
\end{eqnarray*}
\qed

\begin{cor}\label{bim12} We have
	
	
	(i) $\eta_{\lam}^2=\eta_{\lam+\sg}\eta_{\lam-\sg}$, $\sg\in R^0$,$\lam,\sg+\lam\in S$,
	
	(ii) $\eta_{n\lam}=\eta^n_\lam$, $\lam\in S$, $n\in\bbbz$,
	
	(iii) $\eta_{\lam+\sg}=\pm\eta_{\lam}\eta_\sg$, $\lam,\sg, \lam+\sg\in S.$
\end{cor}

\proof
(i) Let $\sg\in R^0$, $\lam,\lam+\sg\in S$. Let $\dot\a\in\rd^\times$ and set $\a=\dot\a+\lam$. Then $\a,\a+\sg\in R^\times$, and so by Lemma \ref{Aug23a}(ii),
$\eta_{\dot\a+\lam}^2=\eta_{\dot\a+\lam+\sg}\eta_{\dot\a+\lam-\sg}.$ 
By inserting $\eta_{\dot\a}\eta_{-\dot\a}$ to both sides, we obtain
$\eta_\lam\eta_\lam=\eta_{\lam+\sg}\eta_{\lam-\sg}$, as required.

(ii) Apply part (i) with $\lam=\sg\in S$, to get $\eta_{2\lam}=\eta^{2\lam}$.
Then by Lemma \ref{bim15}(i), $\eta_{3\lam}=\eta_\lam\eta_{2\lam}=\eta^3_\lam$.
Since $\eta^{-1}_\lam=\eta_{-\lam},$ the result follows inductively. 

(iii) Let $\lam,\sg,\lam+\sg\in S$. Replacing $\lam$ with $\lam+\sg$ in (i), and applying (ii), we get $\eta^2_{\lam+\sg}=\eta_{\lam+2\sg}\eta_\lam=(\eta_\lam\eta_\sg)^2.$
\qed

\begin{pro}\label{bim-15}
	Any finite order Cartan automorphism $\Psi$ for $E$ is diagonal.
\end{pro}

\proof To prove that $\Psi$ is diagonal, we fix $\sg=\sg'+\sg''\in R^0$,   
$\sg',\sg''\in S$. By Lemma \ref{bims30}, we must show that $\eta_\sg^\a=\eta_{\sg}^\b$ for $\a,\b\in R^\times$. Now $\sg'=\tau'+\lam'$ and $\sg''=\tau''+\lam''$ for some
$\tau',\tau''\in S_1$, $\lam',\lam''\in\la S_2\ra$. Let $\a\in R^\times$ with $\a+\sg\in R$. Then
$\a=\dot\a+\sg_s+\lam$, for some $\dot\a\in\rd^\times$, $\sg_s\in S_1$, $\lam\in\la S_2\ra$. Therefore, $\a+\sg=\dot\a+\sg_s+\tau'+\tau''+\lam+\lam'+\lam''\in R$.
This implies that $\sg_s+\tau'+\tau''\in S_1$ and $\lam+\lam'+\lam''\in\la S_2\ra.$
Thus 
\begin{eqnarray*}
	\eta^\a_\sg
	&=&
	\eta_{\dot\a+\sg_s+\tau'+\tau''+\lam+\lam'+\lam''}\eta_{-\dot\a-\sg_s-\lam}\\
		&=&
	\eta_{\dot\a+\sg_s+\tau'+\tau''+\lam+\lam'+\lam''}\eta_{-\dot\a-\sg_s-\tau'-\tau''}\eta_{\dot\a+\sg_s+\tau'+\tau''}\eta_{-\dot\a-\sg_s-\lam}\\
&=&
\eta^{\dot\a+\sg_s+\tau'+\tau''}_{\lam+\lam'+\lam''}\eta_{\dot\a+\sg_s+\tau'+\tau''}\eta_{-\dot\a-\sg_s-\lam}\eta_{\dot\a+\sg_s}\eta_{-\dot\a-\sg_s}\\
&=&
\eta^{\dot\a+\sg_s+\tau'+\tau''}_{\lam+\lam'+\lam''}\eta_{\dot\a+\sg_s+\tau'+\tau''}\eta^{-\dot\a-\sg_s}_{-\lam}\eta^{\dot\a}_{\sg_s}\eta_{-\dot\a}\\
(\hbox{by Lemma \ref{red2}})&=&
\eta^{\dot\theta_s}_{\lam+\lam'+\lam''}\eta^{\dot\theta_s}_{\sg_s+\tau'+\tau''}\eta^{-\dot\theta_s}_{-\lam}\eta^{\dot\theta_s}_{\sg_s}\\
(\hbox{by Lemma \ref{b100}})&=&
\eta_{\lam'}\eta_{\lam''}\eta_{\sg_s+\tau'+\tau''}\eta_{\sg_s}.
\end{eqnarray*}
Thus, we are done if we show that $\eta_{\sg_s+\tau'+\tau''}=\eta_{\sg_s}\eta_{\tau'}\eta_{\tau''}$.

Since $\ind(R)=0$, we have $\ind(S_1)=t$.
Then as we have already mentioned, $S_1$ can be written as $S_1=\uplus_{i=1}^t(\tau_i+2\Lam_1)$, where
$\Lam_1=\sum_{i=1}^t\bbbz\tau_i$, and $\tau_i$'s are distinct modulo $2\Lam_1$.  Now, there exist $1\leq i,j,k\leq t$ such that
$\sg_s=\ep_i\tau_i+2\lam,$ $\tau'=\ep_j\tau_j+2\lam'$, $\tau''=\ep_k\tau_k+2\lam''$, for some $\ep_i,\ep_j,\ep_k\in\{0,1\}$, $\lam,\lam',\lam''\in\Lam_1$. 
Depending on whether $\ep_t$'s are $0$ or $1$, we must consider different possibilities. We examine one of these possibilities, as the other cases can be treated similarly. Let for example $\ep_i=\ep_j=\ep_k=1$. Then since $\sg_s+\tau'+\tau''\equiv\tau_i+\tau_j+\tau_k\in S_1$, modulo $2\Lam_1$, and $\ind(S_1)=t$, we conclude that $|\{i,j,k\}|\leq 2$. Assume for example that $\tau_i=\tau_j\not=\tau_k$. Then, using  Lemma \ref{bim15}(i) and Corollary \ref{bim12}(iii), we have
\begin{eqnarray*}
	\eta_{\sg_s+\tau'+\tau''}
	&=&
	\eta_{2\tau_i}\eta_{\tau_k}\eta_{2\lam}\eta_{2\lam'}\eta_{2\lam''}\\
	&=&
	\eta_{\tau_i}^2\eta_{\tau_k}\eta_{2\lam}\eta_{2\lam'}\eta_{2\lam''}\\
(\hbox{since $\tau_i=\tau_j$})	&=&
	\eta_{\tau_i+2\lam}\eta_{\tau_j+2\lam'}\eta_{\tau_k+2\lam''}\\
	&=&
	\eta_{\sg_s}\eta_{\tau'}\eta_{\tau''},
\end{eqnarray*}
as required. This completes the proof.\qed

\begin{DEF}\label{bim18}
	Let $R$ be of type $B_\ell$ with $\ind(R)=0$.  For $\sg\in R^0$, we define $\eta_\sg:=\eta^\a_\sg$, for some $\a\in R^\times$ {with $\a+\sg\in R^\times$.}
	By Proposition \ref{bim-15}, definition of $\eta_\sg$ is independent of the choice of $\a$.	
\end{DEF}

\begin{lem}\label{bim19}  The following holds:
	
	(i)
	$\eta_{\sg+\sg'}=\eta_{\sg}\eta_{\sg'}$, $\sg,\sg'\in S$.
	
	(ii) For $\sg,\sg'\in S$ and $\sg''\in R^0$, $\eta_{\sg+\sg'+2\sg''}=\eta_{\sg}\eta_{\sg'}\eta_{2\sg''}.$
	
	(iii) For $\sg_1,\sg'_1\in S_1$ and $\lam_2\in\la S_2\ra,$ we have $\eta_{\sg_1+\sg'_1+\lam_2}=
	\eta_{\sg_1}\eta_{\sg'_1}\eta_{\lam_2}$.
\end{lem}

\proof (i) We have $S=S_1\oplus\la S_2\ra$ with $\ind(S_1)=t$. Let
$\sg=\ep_i\tau_i+2\sg_i+\lam,\sg'=\ep_j\tau_j+2\sg_j+\lam'\in S$, where
$1\leq i,j\leq t$,  $\sg_i,\sg_j\in\la S_1\ra$, $\lam,\lam'\in\la S_2\ra$ and $\ep_i,\ep_j\in\{0,1\}$. By definition,
$\eta_{\sg+\sg'}=\eta_{\sg+\sg'}^\a$, for some  $\a\in R^\times$.
Let $\dot\a\in\rds$ and set $\a:=\dot\a-\ep_i\tau_i$. Then $\a+\sg+\sg'=\dot\a+2\sg_i+\lam+\ep_j\tau_j+2\sg_j+\lam'\in R$.
Thus
\begin{eqnarray*}
	\eta_{\sg+\sg'}
	&=&
	\eta^{\dot\a-\ep_i\tau_i}_{\ep_i\tau_i+\ep_j\tau_j+2\sg_i+2\sg_j+\lam+\lam'}\\
		&=&
	\eta_{\dot\a+\ep_j\tau_j+2\sg_i+2\sg_j+\lam+\lam'}\eta_{-\dot\a+\ep_i\tau_i}\\
		&=&
	\eta_{\dot\a+\ep_j\tau_j+2\sg_i+2\sg_j+\lam+\lam'}\eta_{\dot\a}\eta_{-\dot\a}\eta_{-\dot\a+\ep_i\tau_i}\\
		&=&
	\eta^{\dot\a}_{\ep_j\tau_j+2\sg_i+2\sg_j+\lam+\lam'}\eta^{-\dot\a}_{\ep_i\tau_i}\\
	&=&
	\eta_{\ep_j\tau_j+2\sg_i+2\sg_j+\lam+\lam'}\eta_{\ep_i\tau_i}\\
(\hbox{by Lemma \ref{bim15}})		&=&
	\eta_{\ep_j\tau_j}\eta_{2\sg_i}\eta_{2\sg_j}\eta_{\lam}\eta_{\lam'}\eta_{\ep_i\tau_i}.
\end{eqnarray*} 
On the other hand
$$\eta_\sg\eta_{\sg'}=\eta_{\ep_i\tau_i+2\tau_i+\lam}\eta_{\ep_j\tau_j+2\sg_j+\lam'}=\eta_{\ep_i\tau_i}\eta_{2\sg_i}\eta_{\lam}\eta_{\ep_j\tau_j}\eta_{2\sg_j}\eta_{\lam'},$$
so
$\eta_{\sg+\sg'}=\eta_\sg\eta_{\sg'}$.

(ii) Let $\sg,\sg',\sg''$ be as in the statement. Then, we have $2\sg'',\sg'+2\sg''\in S$, and the results easily follows from part (i).

(iii) Let $\sg_1,\sg'_1,\lam_2$ be as in the statement. Then  as $\sg'_1+\lam_2\in S$, we have from part (i) and Lemma \ref{bim15}(ii) that 
$$\eta_{\sg_1+\sg'_1+\lam_2}=\eta_{\sg_1}\eta_{\sg'_1+\lam_2}=\eta_{\sg_1}\eta_{\sg'_2}\eta_{\lam_2}.
$$\qed

\begin{pro}\label{aug29a}
	Let $R$ be of type $B_\ell$ with $\ind(R)=0$. Assume that $\Psi$ is a finite order Cartan automorphism for $E$. Then the corresponding map $\psi:R\rightarrow\bbbk^\times$ is a character.
\end{pro}

\proof By Proposition \ref{bim-15}, $\Psi$ is diagonal. So by Lemma \ref{aug30c}, $\psi$ is a core-character. Therefore, we only need to show that $\eta_{\sg+\sg'}
=\eta_{\sg}\eta_{\sg'}$ for $\sg,\sg',\sg+\sg'\in R^0$.
If $\sg,\sg'\in S$, we are done by Lemma \ref{bim19}(i). Otherwise, using Lemma \ref{bim19}(ii)-(iii), we may assume that 
$\sg=\tau_i+\tau_j$ with $1\leq i\not=j\leq t$, and $\sg'=\ep_k\tau_k+\ep_r\tau_r$, with
$1\leq k\not= r\leq t$, $\ep_k,\ep_r\in\{0,1\}$.
If $\ep_k=\ep_r=0$, then $\sg'=0$ and we are done. If $\ep_k=1$ and $\ep_r=0$, we have
$\sg+\sg'=\tau_i+\tau_j+\tau_k\in R^0$. Since $\ind(S_1)=t$, we get $k=i$ or $k=j$, say $k=j$. Then
\begin{eqnarray*}
	\eta_{\sg+\sg'}&=&
	\eta_{\tau_i+\tau_j+\tau_k}=\eta_{\tau_i+2\tau_j}\\ 
	(\hbox{Lemma \ref{bim19}})
	&=&\eta_{\tau_i}\eta_{2\tau_j}
	=
	\eta_{\tau_i}\eta_{\tau_j}^2=\eta_{\tau_i+\tau_j}\eta_{\tau_j}
	=
	\eta_{\sg}\eta_{\sg'}.
\end{eqnarray*}
Finally, if $\ep_k=\ep_r=1$, then as $i\not=j$, $k\not=r$ and $\sg+\sg'\in S$, we conclude that
$\{i,j\}=\{k,r\}$, say $i=k, j=r$. Then
$$\eta_{\sg+\sg'}=\eta_{2\tau_i+2\tau_j}=\eta_{\tau_i}^2\eta_{\tau_j}^2=\eta^2_{\tau_i+\tau_j}=\eta_{\sg}\eta_{\sg'}.$$
\qed

{To conclude the section, we introduce a few more notations here. Let $\cc$ be the group of Cartan automorphisms for $E$, $\chi$ be the group of characters for $R$, and set $\mathcal{I}:=\{\psi\in\Aut(E)\mid\psi_{|_\hh}=-\id_\hh\}$. Also, we let $\cc^{\hbox{fin}}$, $\chi^{\hbox{fin}}$, and $\mathcal{I}^{\hbox{fin}}$ be the finite order elements of $\cc$, $\chi$, and $\mathcal{I}$, respectively. It is worth noting that $\mathcal{I}^{\hbox{fin}}$ corresponds to the set of Chevalley involutions for $E$. The following theorem summarizes the connections among Chevalley involutions, Cartan automorphisms, and characters.}
	
	\begin{thm}\label{cor021} Let $E$ be an extended affine Lie algebra with root system $R$ such that $\rank\; R>1$ and $\ind(R)=0$. Then 
	
		(i) $\cc^{\hbox{fin}}\cong\chi^{\hbox{fin}}$. 
		
		
		(ii) For any $\tau\in\ii^{\hbox{fin}}$, $\tau\cc^{\hbox{fin}}\sub\ii^{\hbox{fin}}$.
	\end{thm}
	
	\proof Part (i) follows from Propositions \ref{lemaug31a}, \ref{bim11} and \ref{aug29a}. 
	Next let $\tau\in\ii^{\hbox{fin}}$ and let $\Psi\in\cc^{\hbox{fin}}$.
	Let $\a\in R$ and $x_\a\in E_\a$. By part (i),
	$\Psi(x_\a)=\eta_\a x_\a$ for some $\eta_\a\in\bbbk^\times$. Then by Lemma \ref{l1}(i),
	$(\tau \Psi)^2(x_\a)=\tau \Psi(\eta_\a\tau( x_\a))=
	\eta_\a\eta_{-\a}\tau^2(x_\a)=x_\a$. Therefore, $\tau \Psi$ has order $2$ and acts as $-\id$ on $\hh$, so $\tau \Psi\in\ii^{\hbox{fin}}$.
	\qed
	
	\begin{rem}\label{cid1}
		(i) Let $E$ be an extended affine root system of rank $>1$. Assume that $E$ is equipped with a Chevalley involution $\tau$. In \cite{AFI22}, the authors show that the core $E_c$ of $E$ admits a Chevalley system with respect to $\tau$, namely a set $\{x_\a\in E_\a\mid\a\in R^\times\}$ such that for $\a\in R^\times$, $[x_\a,x_{-\a}]=h_\a$ and $\tau(x_\a)=x_{-\a}$. Let $\{\bar x_\a\in E_\a\mid\a\in R^\times\}$ be another Chevalley system with respect to a Chevalley involution, say $\bar\tau$.
		Define $\lam_\a\in\bbbk^\times$ by  $\bar x_\a=\lam_\a x_\a$.
Suppose $\ind(R)=0$ and the Cartan automorphism $\Psi:=\bar\tau\tau$ is of finite order. By Theorem \ref{cor021}, $\Psi$ is diagonal and the associated map $\psi$ is a character. Now let $\d\in R^0$ and
$\a,\a+\d\in R^\times$. Then 
\begin{eqnarray*}
	[x_{-\a},x_{-\d+\a}]&=&\tau[x_\a,x_{\d-\a}]\\
		&=&\bar\tau \Psi[x_\a,x_{\d-\a}]\\
		&=&\psi(\d)\lam_{-\a}\lam_{-\d+\a}[\bar x_{-\a},\bar x_{-\d+\a}]\\
		&=&
		\psi(\d)(\lam_{-\a}\lam_{-\d+\a})^2[x_{-\a},x_{-\d+\a}].
		\end{eqnarray*}
Therefore, $(\lam_\a\lam_{\d-\a})^2=\psi(\d)$. This shows that
for $\d\in R^0$, $\a,\b,\a+\d,\b+\d\in R^\times$, we have
$\lam_\a\lam_{\d-\a}=\pm\lam_\b\lam_{\d-\b},$ see \cite[Lemma 6.1(ii)]{AFI22}. 

{(ii) The content of Section \ref{aug30b} shows that the techniques used in the proofs are heavily reliant on the assumption that
the rank of $E$ is greater than $1$. As seen in \cite{AFI22}, the same situation happens when one deals with constructing Chevalley bases for extended affine Lie algebras. This suggests that a different approach should be considered to investigate the rank $1$ situation, namely the type $A_1$.  The same holds for type $B_\ell$, as our proof relies heavily on the index of the involved semilattices in the structure of the root system. We note that in this case the short roots can be described as (non-isotropic roots of) a disjoint union of $A_1$-type extended affine root systems.} 
	\end{rem}

	\section{\bf Characters vs. root-lattice characters}\setcounter{equation}{0}\label{july20a}
In this section, we discuss the possibility of extending a character for an extended affine root system to a root-lattice character. 
	
		\subsection{\bf Extendability in index zero}\setcounter{equation}{0}
		We begin with the following lemma for extended affine roots systems which resembles and 
		extends a familiar result for finite root systems, see \cite[Corollary 10.2]{Hum72}.
	\begin{lem}\label{partial}
		Let $R$ be an extended affine root system, and $\a_1,\ldots,\a_{n+1}\in R^\times$. Then $w_{\a_1}\cdots w_{\a_n}(\a_{n+1})=\ep_1\a_{i_1}+\cdots+\ep_{m}\a_{i_{m}}$ for some $m\in\bbbz_{>0}$,
		$\a_{i_j}\in\{\a_1,\ldots,\a_{n+1}\}$ and $\ep_i\in\{\pm1\}$, such that $\ep_1\a_{i_1}+\cdots\ep_k\a_{i_k}\in R $ for  $1\leq k\leq m$.
	\end{lem}
	
	\proof Assume $\b=w_{\a_1}\cdots w_{\a_{n}}(\a_{n+1})$. We use induction on $n$.
	Assume first that $n=1$, namely $\b=w_{\a_1}(\a_2)$. If $(\a_1,\a_2)=0$, or $\a_1$ and $\a_2$ are linearly dependent, we are clearly done. Otherwise, there exists $\ep_1,\ep_2\in\{\pm1\}$ such that $\{\ep_1\a_1,\ep_2\a_2\}$ is a base for the rank $2$ finite root system $R\cap (\bbbz\a_1\oplus\bbbz\a_2)$.
	Then the claim is known from finite dimensional theory, see \cite{Hum72}. 
	Assume next that $n>1$. Set $\gamma:=w_{\a_{2}}\cdots w_{\a_n}(\a_{n+1})$. By induction steps,
	$\gamma=\ep_{2}\a_{i_2}+\cdots+\ep_p\a_{i_p}$, where $p\in\bbbz_{>0}$, $\a_{i_j}\in\{\a_2,\ldots,\a_{n+1}\}$, $\ep_{i}\in\{\pm 1\}$ and
	$\ep_2\a_{i_2}+\cdots+\ep_k\a_{i_k}\in R$ {for} $2\leq k\leq p.
	$
	
	Now, $\b=w_{\a_1}(\gamma)$. If $(\gamma,\a_1)=0$, or $\gamma$ and $\a_1$ are linearly dependent we are done. Otherwise, we consider the rank $2$ finite root system
	$R\cap(\bbbz\a_1\oplus\bbbz\gamma)$ and we may assume that, up to a plus-minus sign of $\gamma$, $\{\a_1,\gamma\}$ is a base for this root system.
	Then using the realization of rank $2$ finite root systems, we have
	$$\b=w_{\a_1}(\gamma)\in\{\gamma\pm\a_1,\gamma\pm2\a_1,\gamma\pm3\a_1\}.$$
	Correspondingly, $\b=\gamma\pm\a_1$, $\b=\gamma\pm\a_1\pm\a_1$ or $\b=\gamma\pm\a_1\pm\a_1\pm\a_1$, where $\gamma\pm j\a_1\in R^\times$ if $\gamma\pm(j+1)\a_1\in R^\times$. Thus we are done by induction steps.\qed
	
	We recall from \ref{july20b} that if $R$ is simply laced of rank $>1$ or if $R$ is of type $F_4$ or $G_2$, then
	$\ind(R)=0$.
	
	\begin{lem}\label{diag8}
		Let $R$ be an extended affine root system and $\Phi:R\rightarrow\bbbk^\times$ be a finite-order character for $R$. If
		$ind(R)=0$, then	
		%
		$\Phi$ can be extended to a finite-order root-lattice character.
	\end{lem}	 
	
	\proof By \cite[Proposition 4.41]{Az99}, $R^\times$ contains a basis $\pp$ for the root lattice $\la R\ra$ such that $\w_\pp\pp=R^\times$.
	Consider the group homomorphism $\widehat\Phi:\la R\ra\rightarrow \bbbk^\times$ for $\la R\ra$ given by
	$\widehat\Phi(\a)=\Phi(\a)$ for $\a\in\pp$, see Example \ref{sky4}. Using Lemma \ref{partial}, for $\a\in R^\times$ there exists $\a_1,\ldots,\a_n\in\pp$ and 
	$\ep_1,\ldots,\ep_n\in\{\pm1\}$ such that $\a=\ep_1\a_1+\cdots+\ep_n\a_n$, and
	$\ep_1\a_1+\cdots+\ep_k{\a_k}\in R^\times$ for each $1\leq k\leq n$.
	Then
	\begin{eqnarray*}
		\widehat\Phi(\a)=\widehat\Phi(\ep_1\a_1+\cdots+\ep_n\a_n)&=&
		\prod_{i=1}^n\widehat\Phi(\ep_i\a_i)\\
		&=&
		\prod_{i=1}^n\Phi(\ep_i\a_i)\\
		(\hbox{since $\ep_1\a_1+\ep_2\a_2\in R^\times$})	&=&
		\Phi(\ep_1\a_1+\ep_2\a_2)\prod_{i=3}^n\Phi(\ep_i\a_i)\\
		(\hbox{since $\ep_1\a_1+\ep_2\a_2+\ep_3\a_3\in R^\times$})	&=&
		\Phi(\ep_1\a_1+\ep_2\a_2+\ep_3\a_3)\prod_{i=4}^n\Phi(\ep_i\a_i)\\
		&\vdots&\\
		&=&	\Phi(\ep_1\a_1+\cdots+\ep_n\a_n)=\Phi(\a).\qquad\hbox{\qed}
	\end{eqnarray*}

	\begin{pro}\label{diag5}
		Let $E$ be an extended affine Lie algebra with root system $R$. If
		$\ind(R)=0$ then
		any Cartan automorphism $\Psi$ for $E$ corresponds to a finite-order root-lattice character. In particular, if $R$ is a finite or affine root system, then any character for $R$ extends to a root-lattice character.
	\end{pro}
	
	\proof
	By Propositions \ref{aug29a}, there exists a finite-order character $\psi:R\rightarrow \bbbk^\times$ for $R$ such that $\Psi_{|_{E_\a}}=\psi(\a)\id_{|_{E_\a}}$, $\a\in R$. By Lemma \ref{diag8}, $\psi$ can be extended to a root-lattice character and so we are done.\qed
	
	\begin{rem}\label{diag9}
		We may give an alternative proof for Proposition \ref{diag5} as follows.
		By Proposition \ref{aug29a}, there exists a finite-order character $\psi:R\rightarrow \bbbk^\times$  such that $\Psi_{|_{E_\a}}=\psi(\a)\id_{|_{E_\a}}$, $\a\in R$. By \cite[Lemma 6.2.4]{APT23}, $R^\times$ contains a basis $\pp$ of the root lattice $\la R\ra$ such that the root spaces $E_{\pm\a}$, $\a\in\pp$, generate $E_c$ as a Lie algebra. 
		Consider the group homomorphism $\widehat\psi:\la R\ra\rightarrow \bbbk^\times$ given by
		$\widehat\psi(\a)=\psi(\a)$ for $\a\in\pp$. Let $\widehat\Psi$ be the corresponding Cartan automorphism given in Example \ref{sky4}. Since $E_{\pm\a}$, $\a\in\pp$ generate $E_c$,
		$\Psi$ and $\widehat\Psi$ coincide on the core. Moreover, $\widehat\Psi\Psi=\Psi\widehat\Psi$. Thus by
		Lemma \ref{l2}, $\Psi=\widehat\Psi$, and so ${\widehat\psi}_{|_R}=\psi.$
	\end{rem}

	\subsection{\bf Characters of type {$\mathbf{A_1}$}}
	Lemma \ref{diag8} demonstrates that in order to examine the existence of characters lacking extensions to the root lattice, the type $A_1$ is a suitable candidate. To begin with, we present an extended affine root system of type $A_1$. Consider a one-dimensional real vector space $\vd:=\bbbr\dot\a$, and a $\nu$-dimensional real vector space $\v^0:=\sum_{i=1}^\nu\bbbr\sg_i$. Set $\v:=\vd\oplus\v^0$ and $\Lam=\sum_{i=1}^\nu\bbbz\sg_i$. In what follows, for $\sg,\sg'\in\Lam$, we employ the notation $\sg\equiv\sg'$ to signify that $\sg-\sg'\in 2\Lam$. We define a symmetric bilinear map $\fm$ on $\v$ by 
	$$(\dot\a,\dot\a)=2,\andd(\v^0,\v)=\{0\}.$$
	Let $m$ be a non-negative integer and fix
	\begin{equation}\label{diag13}
		0=\tau_0,\tau_1,\ldots,\tau_m\in\Lam\hbox{ such that for }i\not=j,\;\tau_i\not\equiv\tau_j.
	\end{equation} 
	Set
	\begin{equation}\label{diag14}
		S=\cup_{i=0}^m(\tau_i+2\Lam)\andd R=(S+S)\cup(\pm\dot\a+S).
	\end{equation}
	Then $S$ is a semilattice in $\v^0$ and $R$ is an extended affine root system of type $A_1$ in $\v$ with $R^0=S+S$, moreover, any extended affine root system of type $A_1$ can be described this way, see \pref{bims23} and \cite[Theorem II.2.37]{AABGP97}.
	
	\begin{lem}\label{diag12} Consider the extended affine root system $R=R^0\cup(\pm\dot\a+S)$ of type $A_1$, where $S=\cup_{i=0}^m(\tau_i+2\Lam)$ is given by (\ref{diag13}) and (\ref{diag14}). Assume further that 
		\begin{equation*}\tag{$\star$}
			\tau_{i_1}+\cdots+\tau_{i_k}\not\equiv 0\hbox{ for }3\leq k\leq 6\hbox{ and distinct } 1\leq i_1,\ldots,i_k\leq m.
		\end{equation*} 
		Define $\Phi: R^\times\rightarrow\{\pm1\}$ by
		$$\Phi(\pm\dot\a+2\Lam)=1,\;\Phi(\pm\dot\a+\tau_i+2\Lam)=-1\hbox{ for }1\leq i\leq m,$$
		and $\Phi:R^0\rightarrow\{\pm1\}$ by
		$$\Phi(\sg)=\left\{\begin{array}{ll}
			1&\hbox{if }\sg\in2\Lam\cup((S+S)\setminus S),\\
			-1&S\setminus 2\Lam.
		\end{array}
		\right.
		$$
		Then $\Phi:R\rightarrow\{\pm1\}$ is a character for $R$.	
	\end{lem}
	\proof First we note that the assumption ($\star$) implies that 
	\begin{equation}
		\begin{array}{c}
			\label{diag15}\tau_{i_1}+\tau_{i_2}\not\in S,\;\tau_{i_1}+\tau_{i_2}+\tau_{i_3}\not\in R^0,\;
			\tau_{i_1}+\tau_{i_2}+\tau_{i_3}+\tau_{i_4}\not\in R^0,\vspace{2mm}\\
			\hbox{ if $\tau_{i_j}$'s are distinct}.
		\end{array} 
	\end{equation}
	Now since  $\tau_{i_1}+\tau_{i_2}\not\in S$ for $1\leq i_1\not=i_2\leq m$, we get $R^0=(S+S)=2\Lam\uplus(S\setminus 2\Lam)
	\uplus ((S+S)\setminus S).$ In fact if $\sg\equiv\tau_i$ and $\sg'\equiv\tau_j$ are elements of $S$, then
	$$
	\sg+\sg'\equiv\left\{\begin{array}{ll}
		\tau_i+\tau_j\in (S+S)\setminus S&\hbox{if }1\leq i\not=j\leq m,\\
		2\tau_i\in 2\Lam&\hbox{if }i=j,\\
		\tau_j\in S\setminus 2\Lam&\hbox{if }i=0,\;j\not=0.
	\end{array}
	\right.
	$$
	Therefore $\Phi:R\rightarrow\{\pm1\}$ is a well defined map. From definition of $\Phi$ it follows that $\Phi(\a)=\Phi(\b)$ if $\a\equiv\b$.
We proceed to show that $\Phi$ is a character. In the following, when there is symmetry in the indices, we present the argument for one and leave the others. 	

Let $\a,\b\in R$ with $\a+\b\in R$.
	If $\a,\b\in R^\times$, then there exist $i,j$ and $\ep\in\{\pm1\}$ such that
	$\a\equiv\ep\dot\a+\tau_i$, $\b\equiv-\ep\dot\a+\tau_j$. Then 
	
	\begin{eqnarray*}
		\Phi(\a+\b)=\Phi(\tau_i+\tau_j)&=&\left\{\begin{array}{ll}
			1&\hbox{if }i=j\hbox{ or }1\leq i\not=j\leq m,\\
			-1&i=0,\;j>0.\\
		\end{array}
		\right.\\
		&=&
		\Phi(\a)\Phi(\b).
	\end{eqnarray*}
	If $\a\in R^\times$ and $\b\in R^0$, then $\a\equiv\ep\dot\a+\tau_i$ and
	$\b\equiv\tau_j+\tau_k$ for some $\ep\in\{\pm1\}$ and some $0\leq i,j,k\leq m$.
	If $|\{i,j,k\}|=3$, then as $\a+\b\in R$ and $\tau_s+\tau_t\not\in S$ if $s\not=t$, we get $1\leq i,j,k\leq m$ and
	$\tau_i+\tau_j+\tau_k\in S\setminus 2\Lam$. Thus 
	$\Phi(\a)=-1$, $\Phi(\b)=1$ and $\Phi(\a+\b)=-1$.
	
	Next, assume that $|\{i,j,k\}|=2$. We have the following possibilities, up to symmetry between $j,k$:
	$$j=k,\;i=0\Rightarrow\Phi(\a)=1,\;\Phi(\b)=1,\;\Phi(\a+\b)=1,$$
	$$j=k,\;i\not=0\Rightarrow\Phi(\a)=-1,\;\Phi(\b)=1,\;\Phi(\a+\b)=-1,$$
	$$j\not=k,\;i=j=0\Rightarrow\Phi(\a)=1,\;\Phi(\b)=-1,\;\Phi(\a+\b)=-1,$$
	$$j\not=k,\;i=j\not=0\Rightarrow\Phi(\a)=-1,\;\Phi(\b)=1,\;\Phi(\a+\b)=-1.$$
	If $|\{i,j,k\}|=1$, we have the following possibilities:
	$$i=j=k=0\Rightarrow\Phi(\a)=1,\;\Phi(\b)=1,\;\Phi(\a+\b)=1,$$
	$$i=j=k\not=0\Rightarrow\Phi(\a)=-1,\;\Phi(\b)=1,\;\Phi(\a+\b)=-1.$$
	The above computations show that $\Phi(\a+\b)=\Phi(\a)\Phi(\b)$ for $\a\in R^\times$ $\b\in R^0$, $\a+\b\in R$.
	
	Finally, we consider the case $\a,\b,\a+\b\in R^0=S+S$. We have $\a\equiv\tau_i+\tau_j$, $\b\equiv\tau_k+\tau_t$ for
	some $0\leq i,j,k,t\leq m$. We have the following possibilities, up to symmetry between $i,j$ or $k,t$:
	$$i=j,\;k=t\Rightarrow \Phi(\a)=1,\;\Phi(\b)=1,\;\Phi(\a+\b)=1,$$
	$$i=j=k\not=t=0\Rightarrow \Phi(\a)=1,\Phi(\b)=-1,\;\Phi(\a+\b)=-1,$$
	$$i=j=k=0,\;t\not=0\Rightarrow \Phi(\a)=1,\;\Phi(\b)=-1,\;\Phi(\a+\b)=-1,$$
	$$i=j=0,\;k=t\not=0\Rightarrow \Phi(\a)=1,\;\Phi(\b)=1,\;\Phi(\a+\b)=1,$$
	$$i=j=0,\;k,t\not=0,\;k\not=t\Rightarrow \Phi(\a)=1,\;\Phi(\b)=1,\;\Phi(\a+\b)=1,$$
	$$i=0,\;j=k=t\not=0\Rightarrow \Phi(\a)=-1,\;\Phi(\b)=1,\;\Phi(\a+\b)=-1,$$		
	$$i=0,\;j,k,t\not=0,\;j=k\not=t\Rightarrow \Phi(\a)=-1,\;\Phi(\b)=1,\;\Phi(\a+\b)=-1,$$
	$$i=0,\;j,k,t\not=0,\;|\{j,k,t\}|=3\Rightarrow \Phi(\a)=-1,\;\Phi(\b)=1,\;\a+\b\equiv\tau_j+\tau_k+\tau_t\not\in R^0,$$
	see (\ref{diag15}) for the conclusion that $\a+\b\not\in R^0$ in the last implication.
	The cases we discussed exhaust situations in which  $i=j,\;k=t$ or at least one of $i,j,k,t$ is zero.

	Next we consider the remaining  cases, namely  $1\leq i,j,k,t\leq m$, and $|\{i,j,k,t\}|\geq 2$:
	$$i=j=k\not=t\Rightarrow \Phi(\a)=1,\;\Phi(\b)=1,\;\Phi(\a+\b)=1,$$
	$$i=j, k\not=t,\;\{i,j,k\}|=3\Rightarrow \Phi(\a)=1,\;\Phi(\b)=1,\;\Phi(\a+\b)=1,$$
	$$i\not=j= k\not=t,\;\{i,j,t\}|=3\Rightarrow \Phi(\a)=1,\;\Phi(\b)=1,\;\Phi(\a+\b)=1,$$
	$$|\{i,j,k,t\}|=4\Rightarrow \Phi(\a)=1,\;\Phi(\b)=1,\;\a+\b\equiv\tau_i+\tau_j+\tau_k+\tau_t\not\in R^0,$$
	see (\ref{diag15}) for the conclusion that $\a+\b\not\in R^0$ in the last implication.\qed

	We continue the section by showing that a character for $R$ may fail to extend to a character for $\la R\ra$. 
	\begin{cor}\label{diag10}
		There exist certain characters for extended affine root systems that are not extendable to root-lattice character.
	\end{cor}
	
	\proof Let $R$ be an extended affine root system of type $A_1$ as in  (\ref{diag14}) with $\nu\geq 6$, $\tau_i=\sg_i$ for $1\leq i\leq 6$ and
	$\tau_7=\sg_1+\cdots+\sg_6$. Note that with our choices of $\tau_i$'s, the conditions ($\star$) in the statement of Lemma \ref{diag12} is fulfilled. Now if the character $\Phi$ is extendable to a character $\bar\Phi$ for $\la R\ra$, then we would have $\bar\Phi(\tau_7)=\prod_{i=1}^6\bar\Phi(\sg_i)=\prod_{i=1}^6\Phi(\sg_i)=1$. On the other hand
	by definition of $\Phi$, we have $\bar\Phi(\tau_7)=\Phi(\tau_7)=-1$, a contradiction.\qed
	

\end{document}